\documentclass[12pt]{amsart}
\usepackage{placeins}
\usepackage{amssymb}
\newtheorem{theorem}{Theorem}

\newtheorem{proposition}[theorem]{Proposition}

\newtheorem{corollary}[theorem]{Corollary}

\newtheorem{definition}{Definition}
\newtheorem{lemma}[theorem]{Lemma}

\newcommand{\F}{{\mathcal F}}

\newcommand{\M}{{\mathfrak M}}
\newcommand{\N}{{\mathcal N}}
\newcommand{\Z}{{\mathbb Z}}
\newcommand{\Q}{{\mathbb Q}}
\newcommand{\R}{{\mathbb R}}

\newcommand{\RP}{{\mathbb RP}}

\begin{document}

\unitlength 1.7mm

{

\title{Trinitary algebras}

\author{V.A.~Vassiliev}
\address{Weizmann Institute of Science}
\email{vavassiliev@gmail.com}

\begin{abstract}
A {\em $k$-trinitary algebra} is any subalgebra of the space of smooth functions $f: M \to {\mathbb R}$ that is distinguished in this space by $k$ independent conditions of the form $f(x_i) = f(\tilde x_i) = f(\hat x_i)$, where $x_i, \tilde x_i,$ and $ \hat x_i $ are distinct points in $ M$, $i=1, \dots, k$, or is approximated by such subalgebras. Trinitary algebras naturally arise in the study of {\em discriminant varieties,} that is, the spaces of singular geometric objects, when the property of being singular is formulated in terms of the simultaneous behavior at three distinct points. The simplest singular objects of this kind are the plane curves with triple self-intersections, see \cite{A}, \cite{MD}.

 The spaces of all $k$-trinitary algebras in $C^\infty(M, {\mathbb R})$ are analogous to the spaces of all ideals of finite codimension, which play the same role in the study of discriminants defined in the terms of a single singular point. These spaces are also analogous to the spaces of {\em equilevel algebras} (see \cite{EA}), which arise in the study of discriminants defined by binary singularities. 

We classify the trinitary algebras up to the codimension four in $C^\infty(S^1, {\mathbb R})$, compute the cohomology rings of their varieties and find the Stiefel--Whitney classes of their canonical normal bundles. We also present a series of $(2k-2)$-dimensional cohomology classes of the spaces of trinitary algebras of codimension $2k$ for any natural $k$.

\end{abstract}

\thanks{This work was supported by the Absorption Center in Science of the Ministry of Immigration and Absorption of the State of Israel} 

\maketitle

\section{Introduction}

\subsection{Chord diagrams, triangular diagrams, equilevel algebras, and trinitary algebras}

{\em Chord diagrams} are collections of pairs of distinct points $x_i, \tilde x_i$ on the circle $S^1$. They appear in knot theory because they count the simplest strata of the {\em discriminant set} of maps $S^1 \to \R^3$, that is, the set of maps that are not smooth embeddings. Each such map with $k$ simple self-intersections is characterized by $k$ chords, namely, by the preimages of all self-intersection points. Each $k$-chord diagram defines a subalgebra of codimension $k$ in $C^\infty(S^1, \R)$ consisting of all maps $f$ such that $f(x_i) = f(\tilde x_i)$ for the points $x_i$ and $ \tilde x_i$ of each chord. More complicated points of the discriminant (that is, maps with multiple self-intersection points, self-tangencies, zeros of the derivative, etc) are related to subalgebras
in $C^\infty(S^1, \R)$ that can be approximated by subalgebras defined by chord diagrams. The space of all these subalgebras completes the space of chord diagrams in a manner similar to how the space of ideals of codimension $k$ in $C^\infty(M, \R)$ completes the space of $k$-point configurations in the manifold $M$.
These spaces of subalgebras are studied in \cite{EA}.

Knot invariants {\em of finite type} (see e.g. \cite{bn}, \cite{bl}, \cite{CDM}, \cite{ks}) can be described in the terms of formal linear combinations of chord diagrams that satisfy certain conditions. The simplest condition requires that no chord diagram participating in such a linear combination can include chords whose points $x_i$ and $ \tilde x_i$ are not separated in $S^1 $ by the points of other chords of the same chord diagram. The simplest chord diagram that satisfies this condition is the cross \
{ \unitlength=1.40mm \special{em:linewidth 0.4pt}
\linethickness{1pt}
\begin{picture}(3.33,4.34)
\put(-1.1,3.43){\line(1,-1){3.5}}
\put(-1.1,-0.23){\line(1,1){3.5}}
\thicklines
\put(0.67,1.67){\circle{5.2}}
\end{picture}
} . 
This chord diagram defines the knot invariant equal to the coefficient of the Alexander polynomial at the monomial $t^2$; it is also called the {\em Casson invariant}.

The theory of {\em doodles} is parallel to knot theory: it classifies and investigates 
 the maps $f: S^1 \to \R^2$ having no triple (or more complicated)
points. Unlike {\em generic immersions} $S^1 \to \R^2$ studied in \cite{A}, doodles can have standard self-tangency or simple cusp points because the maps with these singularities are not limits of maps with triple points. However, a doodle $f:S^1 \to \R^2$ cannot have a point $x \in S^1$ at which $f'(x)=f''(x)=0$, nor can it have two points $x$ and $y$ with $f(x)=f(y)$ and $f'(x)=0$. 

\medskip
\noindent
{\bf Remark.} 
1. A slightly different class of plane curves, also called doodles, was considered in \cite{FT}. 
 
2. The theory of {\it immersions} $S^1 \to \R^2$ without triple points is analogous to  framed knot theory or knot theory without the first Reidemeister move.
\medskip

The simplest invariant of doodles, which is fairly analogous to the Casson invariant, was discovered by A.~Merkov, see \cite{Merx}. 

The natural analog of chord diagrams in this theory is {\em triangular diagrams}, that is, collections of triplets of distinct points in $S^1$, that are the preimages of the triple points in $f(S^1)$, see \cite{V-11}. Each triplet $(x, \tilde x, \hat x)$ defines a subspace of codimension two in $C^\infty(S^1, \R^2)$ by the condition $f(x)=f(\tilde x) = f(\hat x)$.
The interpretation of the Merkov's invariant in the terms of triangular diagrams is given in \cite{MD}. Similarly to the case of knots, it is expressed by the simplest triangular diagram \ \
{ \unitlength=1.40mm \special{em:linewidth 0.4pt}
\linethickness{1pt}
\begin{picture}(3.33,4.34)
\put(0.67,-1.00){\line(-3,5){2.33}}
\put(0.67,-1){\line(3,5){2.33}}
\put(-1.66,3.00){\line(1,0){4.67}} 
\put(0.67,4.34){\line(-3,-5){2.33}} 
\put(0.67,4.34){\line(3,-5){2.33}}
\put(-1.67,0.34){\line(1,0){4.67}}
\thicklines
\put(0.67,1.67){\circle{5.2}}
\end{picture}
} ,
that has no triangles, some two vertices of which are not separated in the source circle
by the vertices of other triangles. Continuing the analogy with knot theory, we study the {\em trinitary subalgebras} in $C^\infty(S^1, \R)$, that is, subalgebras that are the limits of algebras of codimension $2k$ defined by $k$ independent triangles. For instance, such an algebra of codimension three can be defined by two points $a, b \in S^1$ and a point $(\alpha: \beta: \gamma) \in \RP^2$, and consist of all functions $f: S^1 \to \R$ such that 
$$f(a)=f(b), \qquad f'(a)=0, \qquad (f'''(a): f''(a): f'(b)) = (\alpha: \beta: \gamma). $$
The space $\overline{TD}_k(S^1)$ 
of all trinitary subalgebras of codimension $2k$ in $C^\infty(S^1, \R)$ has the natural structure of a $3k$-dimensional semialgebraic variety.

\medskip
\noindent
{\bf Example.}
\rm The space $\overline{TD}_1(S^1)$ is homeomorphic to $Sym^3(S^1) \equiv (S^1)^3/S(3)$. In particular, it is homotopy equivalent to $S^1$. 
Its convenient cell decomposition consists of eight cells \ \
{
\unitlength 0.9mm
\begin{picture}(8,7)
\put(4,2){\circle{8}}
\put(8,2){\circle*{1.25}}
\put(0,2){\line(4,3){5.2}}
\put(0,2){\line(4,-3){5.2}}
\put(5.2,5.9){\line(0,-1){7.8}}
\end{picture} \ \ \
\begin{picture}(8,6)
\put(4,2){\circle{8}}
\put(8,2){\circle*{1.25}}
\put(8,2){\line(-4,3){5.2}}
\put(8,2){\line(-4,-3){5.2}}
\put(2.8,5.9){\line(0,-1){7.8}}
\end{picture} \ \ \
\begin{picture}(8,6)
\put(4,2){\circle{8}}
\put(8,2){\circle*{1.2}}
\put(4,6){\line(0,-1){8}}
\put(4,6){\makebox(0,0)[cc]{$\ast$}}
\end{picture} \ \ \
\begin{picture}(8,6)
\put(4,2){\circle{8}}
\put(8,2){\circle*{1.2}}
\put(4,6){\line(0,-1){8}}
\put(4,-2){\makebox(0,0)[cc]{$\ast$}}
\end{picture} \ \ \
\begin{picture}(8,6)
\put(4,2){\circle{8}}
\put(8,2){\circle*{1.2}}
\put(0,2){\line(1,0){8}}
\put(0,2){\makebox(0,0)[cc]{$\ast$}}
\end{picture} \ \ \
\begin{picture}(8,6)
\put(4,2){\circle{8}}
\put(8,2){\circle*{1.2}}
\put(0,2){\line(1,0){8}}
\put(8.8,2){\makebox(0,0)[cc]{$\ast$}}
\end{picture} \ \ \ \
\begin{picture}(8,6)
\put(4,2){\circle{8}}
\put(8,2){\circle*{1.2}}
\put(-0.8,2){\makebox(0,0)[cc]{$\ast$}}
\put(0.8,2){\makebox(0,0)[cc]{$\ast$}}
\end{picture} \ \ \
\begin{picture}(8,6)
\put(4,2){\circle{8}}
\put(8,2){\circle*{1.2}}
\put(7.2,2){\makebox(0,0)[cc]{$\ast$}}
\put(8.8,2){\makebox(0,0)[cc]{$\ast$}}
\end{picture}
} \ of dimensions $3, 2, 2, 2, 1, 1, 1, 0$. 
Here the first picture shows the three-dimensional family of algebras defined by the conditions $f(a)=f(b)=f(c)$ where $a, b,$ and $c$ are all distinct from each other and from the marked point $\bullet$ of the circle. Segments with an asterisk at an endpoint denote the algebras defined by the conditions $f(a)=f(b)$ and $f'(a)=0$. These algebras fill four cells of dimensions two and one, depending on the disposition of the points $a, b,$ and $\bullet$. The double asterisk at point $a \in S^1$ denotes the algebra defined by conditions $f'(a)=f''(a)=0$.
\medskip

We define the structure of a CW-complex on the space $\overline{TD}_2(S^1)$ of all trinitary subalgebras of codimension four and compute its mod 2 cohomology ring.

\begin{theorem}
\label{mainthm1}
The group $H^i(\overline{TD}_2(S^1), \Z_2)$ is isomorphic to $\Z_2$ \ if \ $i=0$ 
or $1$,  to $ \Z_2^2$ \ if $i=2$ or $3$, and is trivial for all $i \geq 4$. 
\end{theorem}

\begin{definition} \rm
The {\it canonical normal bundle} $\N_k$ on the space $\overline{TD}_k(M)$ is the $2k$-dimensional vector bundle whose fiber over any point of this space is the quotient space of $C^\infty(M, \R)$ by the corresponding subalgebra.
\end{definition}

\begin{theorem}
\label{mainthm2}
The group $H^1(\overline{TD}_2(S^1), \Z_2) $ is generated by the first Stiefel--Whitney class $w_1$ of the canonical normal bundle $\N_2$. The group 
$H^2(\overline{TD}_2(S^1), \Z_2) $ is generated by the Stiefel--Whitney classes $w_2$ and $w_1^2$ of this bundle. The classes $w_3(\N_2),$ $w_1(\N_2) \smile w_2(\N_2),$ and $w_1^3(\N_2) \in H^3(\overline{TD}_2(S^1), \Z_2)$ are trivial. 
\end{theorem} 

\begin{corollary}
\label{corrat}
The group $H^i(\overline{TD}_2(S^1), \Q) $ is isomorphic to $\Q$ if $i=0$ or $3$ and is trivial for all other $i$.
\end{corollary}

\noindent
{\it Proof.} The dimensions of rational cohomology groups are not greater than these of $\Z_2$-cohomology groups. By Theorem \ref{mainthm1}, the Euler characteristic of this space is equal to zero. By Theorem \ref{mainthm2}, the Bockstein operator $\mbox{Sq}^1$ acts non-trivially on the generator of the group $H^1(\overline{TD}_2(S^1), \Z_2)$, therefore 
$H^1(\overline{TD}_2(S^1), \Q) \simeq 0$. \hfill $\Box$

\begin{theorem}
\label{mainthm3}
For any $k$, the Stiefel--Whitney class $w_{2k-2}(\N_k) \in H^{2k-2}(\overline{TD}_k(S^1), \Z_2)$ is non-trivial. 
\end{theorem}

Theorem 1 will be proved in \S~\ref{boundaries}, Theorem 2 in \S~\ref{pro2}, Theorem 4 in \S~\ref{series}.

\section{Formal definition and basic properties of trinitary algebras}
\subsection{Main definition}
\label{maindef}

Let $M$ be a smooth real algebraic manifold. Consider the {\em ordered configuration space} $F(M, 3k) \subset M^{3k},$ that is, the set of all sequences $(x_1, \tilde x_1, \hat x_1; x_2, \tilde x_2, \hat x_2; $ $\dots; x_k, \tilde x_k, \hat x_k)$ of $3k$ pairwise distinct points of $M$. Each point of this space defines a subalgebra of codimension $2k$ in $C^\infty(M, {\mathbb R})$ consisting of functions that satisfy all $2k$ equations $f(x_i)=f(\tilde x_i)=f(\hat x_i)$, $i=1, \dots, k$.

Consider an arbitrary algebraic parametric curve $\nu: [0,\varepsilon) \to M^{3k}$ such that all points $\nu(t),$ $t \in (0, \varepsilon),$ belong to the subspace $F(M, 3k) \subset M^{3k}$. For any $t \in (0,\varepsilon)$, denote by $\F_t$ the subalgebra defined by the configuration $\nu(t)$. 

Define $\F_0$ as the space of all functions $f \in C^\infty(M, \R)$ such that there exists a continuous (in $C^{3k}$ topology) family of functions $f_t  \in C^\infty(M, \R)$, $t \in [0,\varepsilon)$, such that $f=f_0$ and $f_t \in \F_t$ for $t \in (0, \varepsilon)$.

 Clearly, the space $\F_0$ is also a subalgebra in $C^\infty(M, {\mathbb R}).$ 

By definition, {\em trinitary algebras} of class $\overline{TD}_k(M)$ are arbitrary subalgebras $\F_0 \subset C^\infty(M, {\mathbb R})$ that can be defined by algebraic curves $\nu$ in this way. 
\medskip

\noindent
{\bf Remark} 
In general, the limit subalgebra $\F_0$ does not uniquely define the point $\nu(0) \in M^{3k}$, even up to the reorderings of $k$ triples $(x_i, \tilde x_i, \hat x_i)$ and of the points within these triples. For instance, the subalgebra of class $\overline{TD}_2(S^1)$ consisting of all functions that satisfy the conditions $f(a)=f(b)=f(c)=f(d)=f(e)$ for five distinct points $a, b, c, d, e \in S^1$ can be defined by this construction with the point $\nu(0) \in (S^1)^6$ equal to $(a, b, c; a, d, e)$ or $(a, b, c; b, d, e)$, or $(a, d, e; b, c, e)$, etc. 
\medskip

 However, the set of all points in $M$ involved in the sequence $\nu(0)$ for any such construction of a trinitary subalgebra is the same for all its constructions; this set is called the {\em support} of the subalgebra. Clearly, the support of a trinitary algebra of codimension $2k$ consists of at most $3k$ points. 

\begin{proposition}
\label{trui}
Suppose that a subalgebra ${\mathcal F}_0 \subset C^\infty(M, {\mathbb R})$ of class $\overline{TD}_k(M)$ has the above construction, and its support consists of the points $a_1, \dots, a_r$, where each point $a_i$ participates in the sequence $\nu(0)$ with some multiplicity $m_i$, \ $m_1 + \dots + m_r =3k$. Then, the subalgebra ${\mathcal F}_0$ contains 

a$)$ all constant functions,

b$)$ the ideal 
\begin{equation} \M_{a_1}^{m_1} \cap \dots \cap \M_{a_r}^{m_r}
\label{ide}
\end{equation} in $C^\infty(M, {\mathbb R})$ consisting of all functions that vanish at each point $a_i$ with multiplicity $m_i$ $($that is, together with all partial derivatives up to degree $m_i-1)$.
\end{proposition}

\noindent
{\it Proof} repeats the proof of Proposition 21 in \cite{EA}. \hfill $\Box$

\subsection{Fighting the infinite-dimensional issues}
\label{findim}

\begin{definition} \rm
A vector subspace $ \F \subset C^\infty(M, \R)$ is {\em of codimension $q$} if, for any vector subspace ${\mathcal L}^n$ of finite dimension $n \geq q$, the intersection $\F \cap {\mathcal L}^n$ is at least $(n-q)$-dimensional, and there exist $q$-dimensional vector subspaces ${\mathcal L}^q$ intersecting $\F$ at the origin only. The $n$-dimensional subspace ${\mathcal L}^n,$ $ n \geq q$, is {\it transversal} to a subspace of codimension $q$ if the dimension of their intersection is exactly $n-q$. 
\end{definition}

Clearly, if the finite-dimensional subspace ${\mathcal L}$ is transversal to a subspace ${\mathcal F}$, then ${\mathcal L}$ is also transversal to any larger vector subspace ${\mathcal F}' \supset {\mathcal F}$. 

Let ${\mathcal L}$ be an arbitrary finite-dimensional vector subspace of $C^\infty(M, {\mathbb R})$ such that 
\begin{equation} 
\label{condition}
\mbox{
\fbox{
\parbox{13cm}{
for any $3k$ different points of $M$ and any collection of $(3k-1)$-jets of functions $M \to {\mathbb R}$ at these points, there exists a function $f \in {\mathcal L}$ whose $(3k-1)$-jets at these $3k$ points coincide with the given ones.}
}
}
\end{equation}

\begin{lemma}
\label{le24}
1. Every trinitary algebra defined as in \S \ref{maindef} by an algebraic curve $\nu: [0, \varepsilon) \to M^{3k}$ is a subspace of codimension $2k$ in $C^\infty(M, \R)$.

2. Each such trinitary algebra intersects the subspace ${\mathcal L}$ along a subspace of codimension $2k$ in ${\mathcal L}$.

3. The intersections of different trinitary algebras with the subspace ${\mathcal L}$ are different subspaces of codimension $2k$ in ${\mathcal L}$.
\end{lemma}

\noindent
{\it Proof} repeats the proof of Lemma 22 in \cite{EA}. \hfill $\Box$
\medskip

Thus, the space $\overline{TD}_k(M)$ can be identified with a semialgebraic subvariety of the finite-dimensional Grassmann manifold $G_{\dim {\mathcal L}-2k}({\mathcal L})$: each subalgebra is associated with its intersection with the space ${\mathcal L}$. If $M$ is compact then this subvariety is also compact. The standard topological structure of the space $\overline{TD}_k(M)$ is induced by this identification from the topology of the Grassmann manifold. This topological structure does not depend on the choice of the space ${\mathcal L}$ satisfying the condition (\ref{condition}). Indeed, each larger space ${\mathcal L}' \supset {\mathcal L}$ defines the same topology of the variety $\overline{TD}_k(M)$, and any two finite-dimensional subspaces are subspaces of a common larger finite-dimensional subspace.

\section{Cell structure of variety $\overline{TD}_2(S^1)$}
\label{cells}

\subsection{} We divide the trinitary subalgebras of $C^\infty(S^1, {\mathbb R})$ into two {\em types}: those algebras
whose supports do not contain or do contain the distinguished point $\bullet = \{0, 2\pi\}\in S^1 \equiv [0,2\pi]/\{0 , 2\pi\}$. We describe 122 families of such subalgebras of codimension four of either of these two types.

All of these families are diffeomorphic to open balls of certain dimensions. They define the structure of a $CW$-complex on the space $\overline{TD}_2(S^1)$. All subalgebras that can be induced from each other via an orientation-preserving diffeomorphism of the pair $(S^1, \bullet)$ to itself belong to the same cell. Thus, the points of their supports (except for $\bullet$) are among the parameters of these cells. Moreover, uniform (in a sense) algebras with the same support also lie in the same cell, and their moduli provide the remaining part of the parameters of the cells.

Each cell of the first type (that is, consisting of algebras whose supports do not contain the point $\bullet$) is related to a cell of the second type. The elements of the latter cell are the algebras obtained from the former algebras by shifting the argument space $S^1$ along its orientation until the support point with the largest coordinate value reaches the distinguished point $\bullet$ . The dimension of the resulting cell is one less than that of the initial cell. The notation of the resulting cell is obtained from the notation of the initial cell by drawing a bar above it: $A_{i j} \mapsto \bar A_{i j}$, $\Omega \mapsto \bar \Omega$, etc.

\subsection{About the pictures}
\label{onnot}

In the figures \ref{A}--\ref{nabla} illustrating the trinitary algebras of first type, the open part $(0, 2\pi)$ of the source circle $S^1$ is shown as a horizontal line. 
 The algebras of the second type corresponding to them can be denoted by the same pictures, in which the rightmost non-generic point is the point $2\pi \equiv 0 \in S^1$.
A tripod with feet at the points $a,$ $b$, and $c$ denotes the condition $f(a)=f(b)=f(c)$. A single arc with endpoints $a$ and $b$ denotes the condition $f(a)=f(b)$. An asterisk $\ast$ at the point $a$ denotes the condition $f'(a)=0$. A double asterisk denotes the condition $f'(a)=f''(a)=0$. A circled asterisk $\circledast$ at the point $a$ denotes the family of pairs of conditions 
$$ f'(a)=0, \quad f'''(a) = \alpha f''(a)$$
depending on the parameter $\alpha \in {\mathbb R}$ . 

An arc connecting the points $a$ and $b$, which are also the legs of some tripod
(and thus satisfy the condition $f(a)=f(b)$),
denotes a family of additional conditions $f'(a)=\alpha f'(b)$, depending on a real parameter $\alpha \neq 0$. This notation is always accompanied by a plus or minus sign that indicates the sign of $\alpha$ for all algebras forming the cell. The notation for more complicated conditions is explained later.

\subsection{Cell structure}
\begin{theorem} 
\label{mainprop}
The variety $\overline{TD}_2(S^1)$ has the structure of a CW-complex with 122 cells of the first type $($that is, those with the support set not containing the distinguished point $\bullet)$ and 122 cells of the second type. 
Namely, it has 
\begin{itemize}
\item[--] ten cells $A_{i j}$ of dimension 6, see Fig.~\ref{A};
\item[--] 
37 cells of dimension 5, including 20 cells $C_{i j}$ $($see Fig.~\ref{C}$)$, one cell $B$ and four cells $E_+$, $E_1,$ $E_2$, and $E_3$ $($see Fig.~\ref{BE}$)$,  and 12 cells $D_{i j}^\pm$ $($see Fig.~\ref{D}$)$;
\item[--] 
45 cells of dimension 4, including four cells $F_i$ $($see Fig.~\ref{F}$)$, 12 cells $G_{i j k l}$ $($see Fig.~\ref{G}$)$, four cells $H_i$ $($see Fig.~\ref{H}$)$, three cells $I_i$ $($see Fig.~\ref{I}$)$, 12 cells $J^{\pm}_{i j l}$ $($see Fig.~\ref{J}$)$, six cells $K_i^{\pm}$ $($see Fig.~\ref{K}$)$, and four cells $W_i^{\pm}$ $($see Fig.~\ref{W}$)$;
\item[--]
23 cells of dimension 3, including six cells $L_{i j l}$ $($see Fig.~\ref{L}$)$, three cells $M_i$ and three cells $N_i$ $($see Fig.~\ref{MN}$)$, two cells $P_i$, two cells $S_i$, two cells $X^{\pm}$ $($see Fig.~\ref{PSX}$)$, four cells $V_i^{\pm}$, and the cell $\Omega$ $($see Fig.~\ref{VOmega}$)$;
\item[--]
six cells of dimension 2, including two cells $Y_i$, two cells $U_i$, one cell $Z$, and one cell $\Theta$ $($see Fig.~\ref{Theta}$)$;
\item[--]
and one cell $\nabla$ of dimension 1, see Fig.~\ref{nabla}.
\end{itemize}
The variety $\overline{TD}_2(S^1)$ also contains the same number of corresponding cells of the second type; their dimensions are 5, 4, 3, 2, 1, and 0, respectively.
\end{theorem}

{\normalsize

\begin{figure}[h]
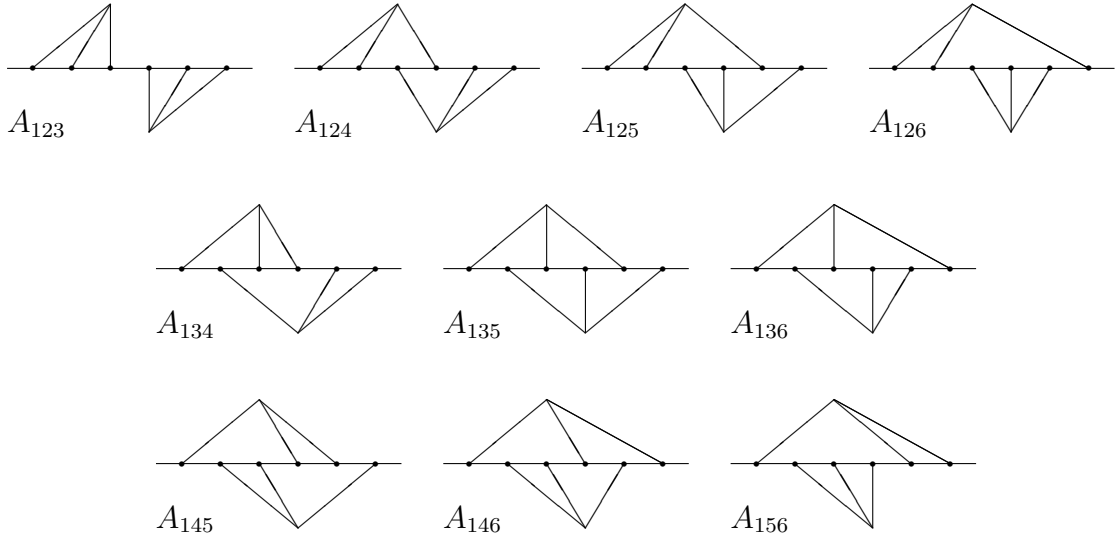

\begin{center}

\end{center}
\caption{Series $A$}
\label{A}
\end{figure}

\begin{figure}[!h]
\begin{center}
%
\caption{$C$}
\label{C}
\end{center}
\end{figure}

\begin{figure}[!h]
\begin{center}
%
\caption{$B, E$}
\label{BE}
\end{center}
\end{figure}

\begin{figure}[h]
%
 
\caption{$D$}
\label{D}
\end{figure}

\unitlength 2mm

\begin{figure}[!h]
%
 
\caption{$F$}
\label{F}
\end{figure}

{\unitlength 1.7mm
\begin{figure}[h]
\begin{center}
%
 
\caption{G}
\label{G}
\end{center}
\end{figure}

\begin{figure}[h]
\begin{center}
%
\caption{$H$} 
\label{H}
\end{center}
\end{figure}

\begin{figure}[h]
\begin{center}
%
 
\caption{$I$}\label{I}
\end{center}
\end{figure}

\begin{figure}[!h]
\begin{center}
%
 \caption{$J$} \label{J}
\end{center}\end{figure}

\begin{figure}[!h]
\begin{center}
%
 \caption{$K$} \label{K} 
\end{center}\end{figure}

\begin{figure}[!h]
\begin{center}
%
 
\caption{W}
\label{W}
\end{center}
\end{figure}

\begin{figure}[!h]
\begin{center}
%
 
\caption{$L$} \label{L} 
\end{center}\end{figure}

\begin{figure}[!h]
\begin{center}
%
 
\caption{$M$ and $N$}
\label{MN}
\end{center}
\end{figure}

\begin{figure}[!h]
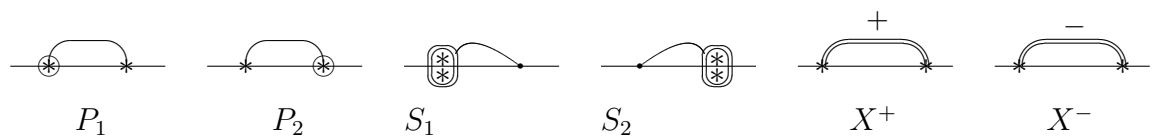

\begin{center}
%
\caption{$P$, $S$, and $X$}
\label{PSX}
\end{center}
\end{figure}

\FloatBarrier

\begin{figure}[!h]
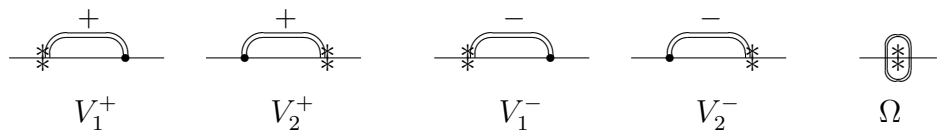

\begin{center}
%
\caption{$V$ and $\Omega$}
\label{VOmega}
\end{center}
\end{figure}

\begin{figure}[!h]
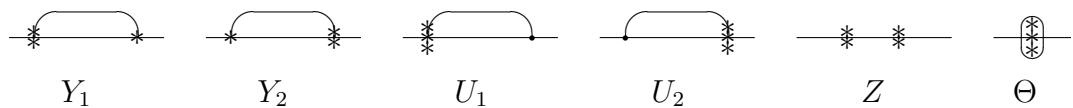

\begin{center}
%
\caption{$Y$, $U$, $Z$, and $\Theta$}
\label{Theta}
\end{center}
\end{figure}

\begin{figure}[h]
\begin{center}
%
\end{center}
\caption{\small $\nabla$}
\label{nabla}
\end{figure}
}
}

The definition of algebras of the classes $A$, $B$, $C$, $D$, $F$, $G$, $H$, $I$, $K$, $L$, $M$, $N$, $P$, $Y$, $U$, $Z$, and $\nabla$ is clear from the pictures. For example, 
the cell $D_{13}^-$ consists of all algebras distinguished by the conditions $f(a)=f(b)=f(c)=f(d)$, $f'(a)= \alpha f'(c),$ where $a<b<c<d \in (0,2\pi)$ and $\alpha <0$.

The algebras of class $E$ depend on three points $a<b<c \in (0, 2\pi)$
and on the parameter $(\alpha: \beta: \gamma) \in \RP^2$,
$\alpha \beta \gamma \neq 0$. Namely, the algebra $E((\alpha:\beta:\gamma), a, b, c)$
consists of functions $f$ such that $$f(a)=f(b)=f(c) \quad \mbox{and} \quad 
(f'(a): f'(b): f'(c)) = (\alpha:\beta:\gamma).$$ The cell
$E_+$ is formed by all algebras of this sort with arbitrary $a<b<c \in (0, 2\pi)$ and the points $(\alpha: \beta:\gamma)$ such that the signs of $\alpha, \beta,$ and $ \gamma$ are the same. The cell $E_1$ (respectively, $E_2$, $E_3$) consists of such algebras with the sign of $\alpha$ (respectively, $\beta, \gamma$) different from the signs of the other two coefficients.

The algebra $J(\alpha; a, b, c)$, $\alpha\neq 0$, of class $J$ consists of functions $f$ satisfying the conditions 
\begin{equation}
\label{eqJ00}
f(a)=f(b)=f(c), \qquad f'(a)=0, \qquad f''(a) = \alpha f'(b).
\end{equation} 
For example, the cell $J_{213}^+$ is the union of such algebras with $0 < b < a < c < 2\pi$ and $\alpha >0$.

The algebras of class $W$ depend on the parameter $(\alpha: \beta: \gamma ) \in \R P^2$, where $\beta \gamma \neq 0$. For any two ordered support points $a$ and $ b$, such an algebra $W((\alpha: \beta: \gamma); a, b)$ consists of all functions $f$ such that $$f(a)=f(b), \qquad f'(a)=0, \qquad ( f'''(a) : f''(a) : f'(b)) = (\alpha : \beta : \gamma). $$ 
These algebras with $a<b$ comprise the cells $W_1^{\pm}$ where the sign $\pm$ is the sign of $\beta/\gamma$. Similarly, the algebras with $a>b$ form the cells $W_2^{\pm}. $ 

The algebras of class $S$ depend on two points $a$, $b$, and the real parameter $\alpha$. The algebra $S(\alpha; a, b)$ consists of functions $f$ such that $$f(a)=f(b) , \qquad f'(a)=f''(a)=0, \qquad f^{\tiny IV}(a)=\alpha f'''(a).$$ These algebras belong to the cell $S_1$ if $a<b$ and to the cell $S_2$ if $a>b$.

The algebras of class $X$ depend on the parameter $\alpha \neq 0$: the algebra $X(\alpha; a, b)$ consists of functions $f$ satisfying the conditions 
$$f(a)=f(b), \qquad f'(a)=0, \qquad f'(b)=0, \qquad f''(a) = \alpha f''(b).$$ Obviously, $X(\alpha; a, b) \equiv X(\alpha^{-1}; b, a)$. These algebras of the first type fill in two cells $X^{\pm}$ depending on the sign of $\alpha$. \label{xx}

The algebras of class $V$ also depend on two points of the circle and one real parameter $\alpha \neq 0$. The algebra $V(\alpha; a, b)$ consists of functions $f$ such that $$f(a)=f(b), \qquad f'(a)=f''(a)=0, \qquad f'''(a) = \alpha f'(b).$$ These algebras of the first type  fill in four cells characterized by the sign of $\alpha$ and by the order of the points $a, b$ in the regular part $(0, 2\pi)$ of the source circle.

The algebras of class $\Omega$ depend on two parameters $\alpha$ and $\beta$. 
The algebra $\Omega (\alpha, \beta; a)$ consists of functions satisfying the conditions 
$$f'(a)=f''(a)=0, \qquad f^{\tiny IV}(a) = \alpha f'''(a), \qquad f^{\tiny V}(a) = \beta f'''(a).$$

The algebras of class $\Theta$ depend on one parameter $\alpha$. The algebra $\Theta(\alpha; a)$ consists of functions $f$ satisfying the conditions $$f'(a)=f''(a)=f'''(a)= 0, \qquad f^{\tiny V}(a) = \alpha f^{\tiny IV}(a).$$ These algebras constitute a two-dimensional cell.

\subsection{The list of algebras of Theorem~\ref{mainprop} is complete}

In \S\S \ref{1sup}--\ref{56sup} we list all possible types of trinitary algebras of codimension four with different cardinalities of their supports. In particular, we show that all of them are counted in Theorem \ref{mainprop}. In \S \ref{boundaries} we present all incidence coefficients (mod 2) of the cells listed in this theorem. It turns out that each of these cells of the first type and of dimension $\leq 5$ occurs with a non-zero coefficient in the formulas for the boundaries of cells of higher dimensions. This implies that all these algebras indeed appear as limits of algebras defined by pairs of triples of points in $S^1$.

\subsubsection{One-point support} 
\label{1sup}
Let $\F$ be a codimension-four trinitary subalgebra defined by a curve germ $\nu: [0,t) \to (S^1)^6$ as in \S \ref{maindef}. 
Suppose that its support is a single point $a$. By Proposition \ref{trui}, $\F \cap \M/ \M_a^6$ is a codimension-four subspace in the space $\M_a/\M_a^6 \simeq \R^5$. 
Each of two families of triples of points participating in its construction converges to $a$ when $t$ tends to 0, therefore $f'(a)=f''(a)=0$ \, for any $f \in \F$. Thus, the line $(\F \cap \M_a)/\M_a^6$ lies in the three-dimensional space $\M^3_a/\M^6_a$. It can be described by a point $(\alpha: \beta:\gamma) \in \RP^2:$ all functions $f \in \F$ satisfy the condition $(f^{V}(a): f^{IV}(a): f'''(a)) = (\alpha:\beta:\gamma)$. 
If $\beta = \gamma=0$ then $\F$ is of class $\nabla$. If $\gamma=0, \beta \neq 0$, then $\F$ is of class $\Theta$. If $\gamma \neq 0$ then $\F$ is of class $\Omega$.

{\normalsize
\unitlength 1.6mm
\begin{figure}
\begin{picture}(10,10)
\put(0,5){\line(1,0){12}}
\put(6,10){\line(-1,-1){5}}
\put(6,10){\line(1,-1){5}}
\bezier{150}(6,10)(9,9)(11,5)
\put(6,0){\line(1,1){5}}
\bezier{150}(6,0)(9,1)(11,5)
\bezier{150}(6,0)(8,4)(11,5)
\put(-0.5,5.7){$a$}
\put(11.7,5.5){$b$}
\end{picture} \qquad \quad
\begin{picture}(10,10)
\put(0,5){\line(1,0){12}}
\put(0,5){\line(1,0){12}}
\put(6,10){\line(-1,-1){5}}
\put(6,10){\line(1,-1){5}}
\bezier{150}(6,10)(3,9)(1,5)
\put(6,0){\line(1,1){5}}
\bezier{150}(6,0)(9,1)(11,5)
\bezier{150}(6,0)(8,4)(11,5)
\put(-0.5,5.7){$a$}
\put(11.7,5.5){$b$}
\end{picture} \qquad \quad
\begin{picture}(10,10)
\put(0,5){\line(1,0){12}}
\put(0,5){\line(1,0){12}}
\put(6,10){\line(-1,-1){5}}
\put(6,10){\line(1,-1){5}}
\bezier{150}(6,10)(9,9)(11,5)
\put(6,0){\line(1,1){5}}
\bezier{150}(6,0)(9,1)(11,5)
\put(6,0){\line(-1,1){5}}
\put(-0.5,5.7){$a$}
\put(11.7,5.5){$b$}
\end{picture} \qquad \quad
\begin{picture}(10,10)
\put(0,5){\line(1,0){12}}
\put(0,5){\line(1,0){12}}
\put(0,5){\line(1,0){12}}
\put(6,10){\line(-1,-1){5}}
\put(6,10){\line(1,-1){5}}
\bezier{150}(6,10)(9,9)(11,5)
\put(6,0){\line(1,1){5}}
\bezier{150}(6,0)(3,1)(1,5)
\put(6,0){\line(-1,1){5}}
\put(-0.5,5.7){$a$}
\put(11.7,5.5){$b$}
\end{picture} \qquad \quad
\begin{picture}(10,10)
\put(0,5){\line(1,0){12}}
\put(1,10){\line(0,-1){5}}
\bezier{150}(1,10)(-1,7.5)(1,5)
\bezier{150}(1,10)(3,7.5)(1,5)
\put(11,0){\line(0,1){5}}
\bezier{150}(11,0)(13,2.5)(11,5)
\bezier{150}(11,0)(9,2.5)(11,5)
\put(0.5,3){$a$}
\put(11,5.5){$b$}
\end{picture}
\caption{Degenerate hypergraphs with two vertices}
\label{lim2}
\end{figure}
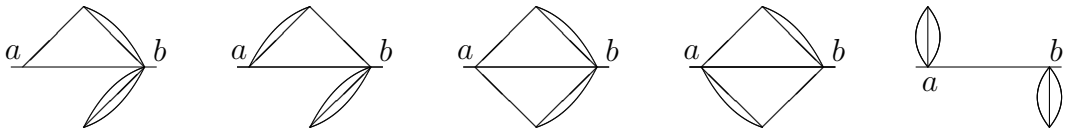
}

\subsubsection{Two-point support} 
\label{2sup}
There are only five combinatorial types of pairs of tripods with two feet to which 
two triples of points constituting the configurations $\nu(t)$ can converge, see Fig.~\ref{lim2}. 

In the case shown in the leftmost picture, we have three conditions $f(a)=f(b),$ $f'(b)=f''(b)=0$. The space $(\F \cap \M_b)/\M^5_b$ is then a subspace of codimension one in the space $\M^3_b/\M^5_b \simeq \R^2$. Thus, the algebra $\F$ is determined by a point $(\alpha:\beta) \in \RP^1$: all functions $f \in \F$ satisfy an additional condition $(f'''(b):f^{IV}(b)) = (\alpha: \beta)$. If $\alpha=0$ then $\F$ is of class $U$, if $\alpha \neq 0$ then $\F$ is of class $S$. 

The second picture in Fig.~\ref{lim2} implies four conditions $f(a)=f(b), $ $f'(a)=0=f'(b)=f''(b)$. Thus, $\F$ is determined by two points $a$ and $b$ and is of class $Y$.

The third picture of Fig.~\ref{lim2} implies two conditions $f(a)=f(b)$ and $f'(b)=0$. The space $(\F \cap \M_a \cap \M^2_b)/(\M^2_a \cap \M^4_b)$ is then a subspace of codimension two in the space $(\M_a \cap \M^2_b)/(\M^2_a \cap \M^4_b) \simeq \R^3$. It is characterized by a point $(\alpha:\beta:\gamma) \in \RP^2$: all functions $f \in \F$ satisfy the additional condition $(f'(a):f''(b):f'''(b))= (\alpha:\beta:\gamma)$. If $\alpha=\beta=0$, then $\F$ is of class $Y$. If $\alpha=0, \beta \neq 0$, then $\F$ is of class $P$. If $\alpha \neq 0, \beta = 0$, then $\F$ is of class $V$. If $\alpha \neq 0 \neq \beta,$ then $\F$ is of class $W$.

The fourth picture of Fig.~\ref{lim2} implies three conditions $f(a)=f(b)$, $f'(a)=0=f'(b)$. The algebra $\F$ is then characterized by the 
space $(\F \cap \M^2_a \cap \M^2_b)/(\M^3_a \cap \M^3_b)$,
which is a line in the plane $(\M^2_a \cap \M^2_b)/(\M^3_a \cap \M^3_b)$. Thus, $\F$ is defined by a point $(\alpha:\beta) \in \RP^1$: all functions $f \in F$ satisfy the fourth condition $(f''(a):f''(b))= (\alpha:\beta)$. If $\alpha=0$ or $\beta=0$ then $\F$ is of class $Y$, if $\alpha \neq 0 \neq \beta$ then $\F$ is of class $X$.

The last picture of Fig.~\ref{lim2} implies four conditions $f'(a)=f''(a)=f'(b)=f''(b)=0$, so algebra $\F$ is determined by the points $a$ and $b$ and is of class $Z$. 

{\normalsize
\unitlength 1.6mm
\begin{figure}
\begin{picture}(10,10)
\put(0,5){\line(1,0){10}}
\put(5,1){\line(1,1){4}}
\put(5,1){\line(0,1){4}}
\put(5,1){\line(-1,1){4}}
\put(5,9){\line(1,-1){4}}
\put(5,9){\line(0,-1){4}}
\put(5,9){\line(-1,-1){4}}
\put(0.5,4){\makebox(0,0)[cc]{\footnotesize $a$}}
\put(4.3,4){\makebox(0,0)[cc]{\footnotesize $b$}}
\put(9.5,4){\makebox(0,0)[cc]{\footnotesize $c$}}
\end{picture} \quad
\begin{picture}(10,10)
\put(0,5){\line(1,0){10}}
\put(5,1){\line(1,1){4}}
\put(5,1){\line(0,1){4}}
\bezier{150}(5,1)(8,2)(9,5)
\put(1,9){\line(0,-1){4}}
\bezier{150}(1,9)(3,7)(1,5)
\bezier{150}(1,9)(-1,7)(1,5)
\put(1,4){\makebox(0,0)[cc]{\footnotesize $a$}}
\put(5,6){\makebox(0,0)[cc]{\footnotesize $b$}}
\put(9,6){\makebox(0,0)[cc]{\footnotesize $c$}}
\end{picture} \quad
\begin{picture}(10,10)
\put(0,5){\line(1,0){10}}
\put(5,1){\line(1,1){4}}
\put(5,1){\line(0,1){4}}
\bezier{150}(5,1)(8,2)(9,5)
\put(5,9){\line(0,-1){4}}
\put(5,9){\line(-1,-1){4}}
\bezier{150}(5,9)(2,8)(1,5)
\put(1,4){\makebox(0,0)[cc]{\footnotesize $a$}}
\put(4.3,4){\makebox(0,0)[cc]{\footnotesize $b$}}
\put(9,6){\makebox(0,0)[cc]{\footnotesize $c$}}
\end{picture} \quad
\begin{picture}(10,10)
\put(0,5){\line(1,0){10}}
\put(5,1){\line(1,1){4}}
\put(5,1){\line(-1,1){4}}
\bezier{150}(5,1)(8,2)(9,5)
\put(5,9){\line(0,-1){4}}
\put(5,9){\line(-1,-1){4}}
\bezier{150}(5,9)(2,8)(1,5)
\put(0.5,4){\makebox(0,0)[cc]{\footnotesize $a$}}
\put(5,4){\makebox(0,0)[cc]{\footnotesize $b$}}
\put(9,6){\makebox(0,0)[cc]{\footnotesize $c$}}
\end{picture} \quad
\begin{picture}(10,10)
\put(0,5){\line(1,0){10}}
\put(5,9){\line(-1,-1){4}}
\bezier{150}(5,9)(6,7)(5,5)
\bezier{150}(5,9)(4,7)(5,5)
\put(5,1){\line(1,1){4}}
\bezier{150}(5,1)(6,3)(5,5)
\bezier{150}(5,1)(4,3)(5,5)
\put(1,4){\makebox(0,0)[cc]{\footnotesize $a$}}
\put(4,4){\makebox(0,0)[cc]{\footnotesize $b$}}
\put(9,6){\makebox(0,0)[cc]{\footnotesize $c$}}
\end{picture}
\quad
\begin{picture}(10,10)
\put(0,5){\line(1,0){10}}
\put(5,1){\line(1,1){4}}
\put(5,1){\line(0,1){4}}
\bezier{150}(5,1)(8,2)(9,5)
\put(5,9){\line(1,-1){4}}
\put(5,9){\line(0,-1){4}}
\put(5,9){\line(-1,-1){4}}
\put(1,4){\makebox(0,0)[cc]{\footnotesize $a$}}
\put(4.3,4){\makebox(0,0)[cc]{\footnotesize $b$}}
\put(9.5,6){\makebox(0,0)[cc]{\footnotesize $c$}}
\end{picture} 
\quad
\begin{picture}(10,10)
\put(0,5){\line(1,0){10}}
\put(5,1){\line(1,1){4}}
\put(5,1){\line(0,1){4}}
\put(5,1){\line(-1,1){4}}
\put(1,9){\line(0,-1){4}}
\bezier{150}(1,9)(3,7)(1,5)
\bezier{150}(1,9)(-1,7)(1,5)
\put(0.7,4){\makebox(0,0)[cc]{\footnotesize $a$}}
\put(5,6){\makebox(0,0)[cc]{\footnotesize $b$}}
\put(9,6){\makebox(0,0)[cc]{\footnotesize $c$}}
\end{picture} 
\caption{Limit hypergraphs with three vertices}
\label{lim3}
\end{figure}
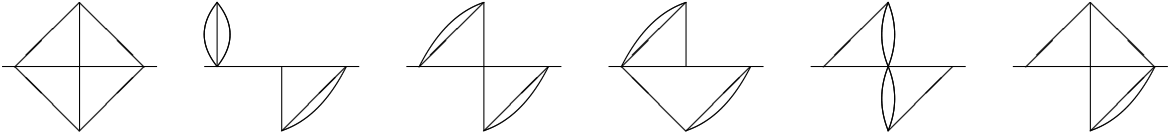
}

\subsubsection{Three-point support} 
\label{3sup}
If the points of the configurations $\nu(t)$ tend to three points $a, b, c$ of the line, then the limit pair of tripods can be of one of the seven types shown in Fig.~\ref{lim3}. 

The first picture of this figure implies two conditions $f(a)=f(b)=f(c)$. The algebra $\F$ is then characterized by the 
space $(\F \cap \M_a \cap \M_b \cap \M_c)/(\M^2_a \cap \M^2_b \cap \M^2_c)$, which is a subspace of codimension two in the space $(\M_a \cap \M_b \cap \M_c)/(\M^2_a \cap \M^2_b \cap \M^2_c) \simeq \R^3$. It is therefore defined by a point $(\alpha: \beta:\gamma) \in \RP^2$: all functions $f \in \F$ satisfy the condition $(f'(a):f'(b):f'(c)) = (\alpha:\beta:\gamma)$. If some two of the coefficients $\alpha, \beta, \gamma$ are equal to zero, then $\F$ is of class $N$. If only one of them is equal to zero, then $\F$ is of class $K$. If none of them are equal to 0 then $\F$ is of class $E$.

The second picture of Fig.~\ref{lim3} implies four conditions $f'(a)=f''(a)=0$, $f(b)=f(c)$, $f'(c)=0$. The algebra $\F$ is determined by the ordered set of points $a, b, c$, and is of class $L$. 

The third and fourth pictures imply four conditions $f(a)=f(b)=f(c)$, $f'(a)=0$, $f'(c)=0$, so the algebra $\F$ is determined by its support points and is of class $N$.

The fifth picture implies three conditions $f(a)=f(b)=f(c)$, $f'(b)=0$. The algebra $\F$ is determined by the image of the subspace $\F \cap \M^2_b$ in the space $\M^2_b/\M^4_b \simeq \R^2$. Let $(\alpha: \beta) \in \RP^1$ be the parameter of this image. The algebra $\F$ then consists of the functions $f $ satisfying additionally the condition $(f''(b):f'''(b)) = (\alpha:\beta)$. If $\alpha=0$ then $\F$ is of class $M$, if $\alpha \neq 0$ then $\F$ is of class $I$. 

The sixth picture implies three conditions $f(a)=f(b)=f(c)$, $f'(c)=0$. The algebra $\F$ is then determined by the class of the line $\F \cap \M_b \cap \M^2_c $ in the space $ (\M_b \cap \M^2_c)/(\M^2_b \cap \M^3_c)$, that is, by the common ratio $(f'(b):f''(c)) \in \RP^1$ for all $f \in \F$. If $f'(b)\equiv 0$, then $\F$ is of class $N$. If $f''(c)\equiv 0$, then it is of class $M$. If neither of these conditions holds for all $f \in \F$, then $\F$ is of class $J$. 

The seventh picture implies four conditions $f(a)=f(b)=f(c)$, $f'(a)=f''(a)=0$ and determines algebras of class $M$.

{\normalsize
\begin{figure}
\begin{picture}(10,12)
\put(0,5){\line(1,0){12}}
\put(6,0){\line(-2,5){2}}
\put(6,0){\line(2,5){2}}
\put(6,0){\line(1,1){5}}
\put(1,10){\line(0,-1){5}}
\bezier{150}(1,10)(3,7.5)(1,5)
\bezier{150}(1,10)(-1,7.5)(1,5)
\put(0.5,3){$a$}
\put(3.5,5.7){$b$}
\put(7.3,5.7){$c$}
\put(10.5,5.7){$d$}
\end{picture} \qquad \quad
\begin{picture}(10,12)
\put(0,5){\line(1,0){12}}
\put(6,0){\line(1,1){5}}
\put(6,0){\line(2,5){2}}
\bezier{150}(6,0)(9,1)(11,5)
\put(6,10){\line(-2,-5){2}}
\put(6,10){\line(-1,-1){5}}
\bezier{150}(6,10)(3,9)(1,5)
\put(0.5,3){$a$}
\put(3.5,3){$b$}
\put(7,5.7){$c$}
\put(10,5.7){$d$}
\end{picture} \qquad \quad
\begin{picture}(10,12)
\put(0,5){\line(1,0){12}}
\put(6,0){\line(1,1){5}}
\put(6,0){\line(2,5){2}}
\put(6,0){\line(-1,1){5}}
\put(6,10){\line(-2,-5){2}}
\put(6,10){\line(-1,-1){5}}
\bezier{150}(6,10)(3,9)(1,5)
\put(0.5,3){$a$}
\put(3.5,3){$b$}
\put(7,5.7){$c$}
\put(10,5.7){$d$}
\end{picture} \qquad \quad 
\begin{picture}(10,12)
\put(0,5){\line(1,0){12}}
\put(6,0){\line(1,1){5}}
\put(6,0){\line(2,5){2}}
\put(6,0){\line(-2,5){2}}
\put(6,10){\line(-2,-5){2}}
\put(6,10){\line(-1,-1){5}}
\bezier{150}(6,10)(3,9)(1,5)
\put(0.5,3){$a$}
\put(3,3){$b$}
\put(7,5.7){$c$}
\put(10,5.7){$d$}
\end{picture} \qquad \quad 
\begin{picture}(10,12)
\put(0,5){\line(1,0){12}}
\put(6,0){\line(1,1){5}}
\put(6,0){\line(2,5){2}}
\put(6,0){\line(-2,5){2}}
\put(6,10){\line(-2,-5){2}}
\put(6,10){\line(-1,-1){5}}
\put(6,10){\line(2,-5){2}}
\put(0.5,3){$a$}
\put(3,3){$b$}
\put(8,5.5){$c$}
\put(10.5,5.5){$d$}
\end{picture} 
\caption{Limit hypergraphs with four vertices}
\label{lim4}
\end{figure}
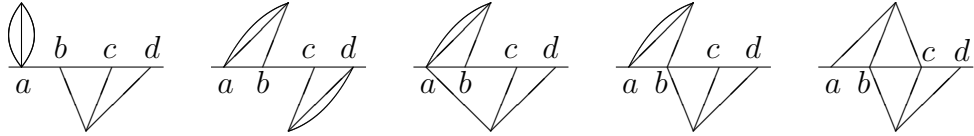
}

\subsubsection{Four-point support}
\label{4sup}

A limit pair of tripods with feet at four points can have only one of the five isomorphism classes shown in Fig.~\ref{lim4}. 

The first picture implies four conditions $f'(a)=f''(a)=0$, $f(b)=f(c)=f(d)$. Thus, the algebra $\F$ is determined by the points of the support and is of class $F$. 

The second picture implies four conditions $f(a)=f(b)$, $f(c)=f(d)$, $f'(a)=0=f'(d)$. The algebra $\F$ then is of class $G$.

The third picture implies four conditions $f(a)=f(b)=f(c)=f(d)$, $f'(a)=0$, and defines an algebra $\F$ of class $H$. The same is true also for the fourth picture.

The fifth picture implies three conditions $f(a)=f(b)=f(c)=f(d)$. Any trinitary algebra $\F$ with this picture satisfies an additional condition of the form $(f'(b):f'(c))= (\alpha:\beta) \in \RP^1$. If one of the coefficients $\alpha$ or $\beta$ is equal to 0, then this algebra is of class $H$; if $\alpha \neq 0 \neq \beta$ then it is of class $D$. 

\subsubsection{Five- and six-point supports} 
\label{56sup}
These cases are trivial and yield algebras of types $B$, $C$, and $A$.

\FloatBarrier

\subsection{End of the proof of Theorem \ref{mainprop}}
The definitions of all cells mentioned in this theorem immediately imply that, for each $i$-dimensional cell $\beta^i_\alpha$, there is a parameterization \begin{equation}
\varphi^i_\alpha: D^i \to \overline{TD}_2(S^1) \subset G_{\dim {\mathcal L}-4}(\mathcal L),
\label{paramcell}
\end{equation}
 where $D^i$ is the standard open unit ball in ${\mathbb R}^i$, and $\varphi_\alpha^i $ is a regular algebraic map that diffeomorphically sends $D^i$ onto $\beta^i_\alpha$. 

\begin{proposition}
The space $\overline{TD}_2(S^1)$ has the structure of a CW-complex with these cells. 
\end{proposition}

\noindent
{\it Proof.} The proof repeats that of Lemma 25 in \cite{EA} and concludes the proof of Theorem \ref{mainprop}. \hfill $\Box$ $\Box$

\subsection{Important filtration}

Let $M$ be an one-dimensional manifold and $\F \in \overline{TD}_k(M)$ be a trinitary subalgebra.

\begin{definition} \rm
The {\it multiplicity} of a support point \ $a$ \ of the algebra $\F$ is the smallest natural number $m$ such that the algebra $\F$ contains all functions in $\M_a^m$ that are identically zero in some neighborhoods of all other support points of $\F$. The {\it multiplicity of the algebra} $\F$ is the sum of the multiplicities of all its support points.
\end{definition}

\begin{proposition}
For any natural number $p$, the set of all trinitary algebras $\F$ of multiplicity $\leq p$ is closed in $\overline{TD}_k(M)$.
\end{proposition}

\noindent
{\it Proof} of this statement repeats the proof of the analogous statement (Proposition 26) on the equilevel algebras, see \cite{EA}.

\begin{proposition}
The algebras of classes $B$, $H$, $M$, $N$, $Y$, $U$, and $\nabla$ are of multiplicity 5, all other algebras from Theorem \ref{mainprop} are of multiplicity 6. \hfill $\Box$
\end{proposition}

An important property of this filtration is that the set of all trinitary algebras of a fixed multiplicity in $\overline{TD}_k(S^1)$ is fibered over the source circle of our functions.
Indeed, we can consider this circle as the Lie group of complex numbers of unit norm. Associate with each trinitary subalgebra $\F \in \overline{TD}_k(S^1)$ the product of its support points taken with their multiplicities.
The resulting map $\overline{TD}_k(S^1) \to S^1$ is continuous on the set of trinitary algebras of any fixed multiplicity $m$. The restriction of this map to this set commutes with the simultaneous rotation of the source and the target circles by the angles \ $\alpha$ \ and \ $m \alpha$, respectively.

\section{A series of homology classes of spaces $\overline{TD}_k(S^1)$ for arbitrary $k$}
\label{series}

Let $a, b, $ and $c$ be arbitrary three points in $S^1$; let $\Delta: (S^1, a) \to (S^1, b)$ and $\tilde \Delta: (S^1, a) \to (S^1, c)$ be two germs of local diffeomorphisms between the neighborhoods of these three points. The algebra $E(\Delta, \tilde \Delta; a, b, c)$ is the space of all functions $f: S^1 \to \R$ such that the $(k-1)$-jets of three functions $f$, $f \circ \Delta$, and $f \circ \tilde \Delta$ at the point \ $a$ \ are the same.

It is easy to see that this space is a subalgebra that depends only on $(k-1)$-jets of the diffeomorphisms $\Delta$ and $\tilde \Delta$. Different pairs of these jets define different algebras. 

\begin{proposition}
Each algebra $E(\Delta, \tilde \Delta; a, b, c)$ belongs to the space $\overline{TD}_{2k}(S^1).$
\end{proposition}

\noindent
{\it Proof.} Take the family of triples $(a, b, c; a+t, \Delta(a+t), \tilde \Delta(a+t); a+2 t, \Delta(a+2 t), \tilde \Delta(a+2 t); \dots; a+(k-1) t, \Delta(a+(k-1) t), \tilde \Delta(a+(k-1) t)) ) \in F(S^1, 3k)$. \hfill $\Box$ \medskip

By Proposition \ref{trui}, any such subalgebra $\F$ is characterized by the image of the space $\F \cap \M_a$ in the quotient space $$\M_a/ (\M^k_a \cap \M^k_b \cap \M^k_c) \simeq \R^{3k-1}.$$ It is a subspace of codimension $2k$ and hence of dimension $k-1$ there. By definition, it lies in the subspace 
$$(\M_a \cap \M_b \cap \M_c)/ (\M^k_a \cap \M^k_b \cap \M^k_c) \simeq \R^{3k-3}$$ and defines a point of the Grassmann manifold $G_{k-1}(\R^{3k-3})$. The set of all such points defined by arbitrary local diffeomorphisms $\Delta$ and $\tilde \Delta$ is a semialgebraic $(2k-2)$-dimensional subvariety in this manifold. Denote by $\Xi(a,b,c)$ the algebraic closure of this space.
Then we have the following detalization of Theorem \ref{mainthm3}.

\begin{proposition}
For any three distinct points $a, b, c \in S^1,$ the Stiefel--Whitney class $w_{2k-2}(\N_k)$ takes the non-zero value on the cycle $\Xi(a,b,c)$.
\end{proposition}

\noindent
{\it Proof.} Let $\{g_1, g_2, g_3\}$ be a generic triple of  functions on $S^1$ which belong to the space ${\mathcal L}$ (see \S~\ref{findim}). Their cosets modulo the trinitary algebras define three sections of the vector bundle $\N_k$ on $\overline{TD}_k(S^1)$. The value of the class $w_{2k-2}(N_k)$ on a $(2k-2)$-dimensional cycle in $\overline{TD}_k(S^1) \subset G_{\dim {\mathcal L} -2k}({\mathcal L})$ is equal to the intersection index of that cycle with the set of points of $G_{\dim {\mathcal L} -2k}({\mathcal L})$ such that the corresponding fiber of the tautological bundle has a positive-dimensional intersection with the three-dimensional space spanned by the functions $g_i$. Tautologically, the
intersection set of these two cycles in $G_{\dim {\mathcal L} -2k}({\mathcal L})$
consists of trinitary algebras of codimension $2k$ that contain a nonzero linear combination ${\bf g}$ of the functions $g_i$.

This linear combination must have equal values at the points $a, b, $ and $c$. Since the system of functions $g_i$ is generic, this condition fixes the linear combination ${\bf g}$ up to a nonzero scalar, and the derivatives of this linear combination at the points $a, b, c$ are not equal to zero. Therefore, the $(k-1)$-jets of local diffeomorphisms $\Delta: (S^1, a) \to (S^1,b) $ and $\tilde \Delta: (S^1, a) \to (S^1,c)$ are uniquely defined by the condition that this linear combination belongs to $E(\Delta, \tilde \Delta; a, b, c)$, that is, ${\bf g} \equiv {\bf g} \circ \Delta \equiv {\bf g} \circ \tilde \Delta (\mbox{mod } \M_a^k)$ in a neighborhood of the point $a$. The algebra $E(\Delta, \tilde \Delta; a, b, c)$ is the unique (and transversal) intersection point of the cycle $\Xi(a,b,c) \subset G_{\dim {\mathcal L} -2k}({\mathcal L})$ and the set of all subspaces of codimension $2k$ in ${\mathcal L}$ not in general position with the subspace spanned by the functions $g_1, g_2,$ and $ g_3$. \hfill $\Box$ $\Box$

{\normalsize
\section{Boundary operators}
\label{boundaries}

Consider the cellular chain complex with coefficients in ${\mathbb Z}_2$ defined by our CW-structure of the variety $\overline{TD}_2(S^1)$. 
Most of the incidence coefficients in the formulas for its boundary operators $C_i \to C_{i-1}$ are obvious. In particular, this is the case if the incidence of cells is induced by the collision of only two simple points of the supports of the subalgebras from the larger cell (see for example the degeneration of cells of type $A_{i j}$ into cells $B$ or $C_{k l}$), or it is induced by the boundary values of the numerical conditions on the derivatives at the points of a fixed configuration (as in the case of the cells $D^\pm_{i j}$ approaching the cells $H_i$). Here is an example of a more complicated collision.

\subsection{Incidence coefficient $[\bar J^{\pm}_{312}, \bar S_2]$.}
\label{incJS}
A typical function of an algebra $J(\alpha; 0, a, b)$ of class $J^-_{312}$, whose  support point \ $a$ \ approaches the critical point at \ 0, \ looks as shown in Fig.~\ref{incid2}.

\begin{figure}
\unitlength 1mm
\begin{picture}(100,20)
\put(0,10){\line(1,0){100}}
\bezier{400}(40,0)(55,17)(70,5)
\bezier{400}(70,5)(85,-5)(100,20)
\put(0,20){\line(1,-1){20}}
\put(57,7.5){\makebox(0,0)[cc]{\footnotesize $0$}}
\put(93,8){\makebox(0,0)[cc]{\footnotesize $a$}}
\put(8,8){\makebox(0,0)[cc]{\footnotesize $b$}}
\end{picture}
\caption{Incidence coefficient $[\bar J^-_{312}, \bar S_2]$}
\label{incid2}
\end{figure}

All functions of such an algebra satisfy the conditions
\begin{equation}
\label{eqJJ} 
 f(a)=f(b)=f(0), \qquad f'(0)=0, \qquad f'(a) = \lambda f''(0),
\end{equation}
where $\lambda = \alpha^{-1}< 0,$ cf. (\ref{eqJ00}). 
Take the coordinate of the point $a$ as the parameter of the collision. Let the other parameter $\lambda$ of the algebras (\ref{eqJJ}) depend on $a$ as 
\begin{equation}
\label{lamed}
 \lambda (a) = p a + q a^2 + \dots.
\end{equation}
 At the end of this collision, we need to obtain a function of the algebra $S(\theta; 0, b),$ all whose elements $f$ satisfy the conditions 
\begin{equation}
\label{eqfin}
 f(b)=f(0), \qquad f'( 0)=f''(0)=0, \qquad    f^{\tiny \mbox{IV}}(0) = \theta f'''(0).
\end{equation}
 with an arbitrary prescribed value of the parameter $\theta$.

 According to Proposition \ref{trui}, this limit algebra is completely characterized by its intersection with the space ${\mathbb R}^3$ of polynomials 
\begin{equation}
A x^2 + B x^3 + C x^4 , 
\label{trans}
\end{equation}
which is transversal to $\M^5_0 $ in the space $\M^2_0$.
Each algebra (\ref{eqJJ}) corresponding to some value of the parameter $a$ intersects the space of polynomials (\ref{trans}) along the set of such polynomials
which satisfy two linear conditions $f(a)=f(0)$ and $f'(a)=\lambda f''(0)$, i.e.,
\begin{equation}
\label{ee1}
\begin{split}
A + B a + C a^2 & =0, \\
2A a + 3 B a^2 + 4 C a^3 & = 2A (p a + q a^2 + \dots) .
\end{split}
\end{equation}
The two planes in ${\mathbb R}^3$ defined by these conditions tend to the plane $\{A=0\}$ when $a$ tends to $0$. The limit position of their intersections, which is a line in the latter plane, satisfies the equations (\ref{eqfin}) with some given value $\theta$ if it consists of functions (\ref{trans}) with $A=0, \frac{C}{B} = \frac{\theta}{4}$.

To compute this limit position, we express the variable \ $A$ \ from one of the equations (\ref{ee1}) and substitute it into the other: 
\begin{equation}
\label{linbcd}
3 B a + 4 C a^2 \equiv -2(B + C a)((p-1) a + q a^2 + ...)
\end{equation}
This is a linear relation between the variables $B$ and $C$, whose coefficients are some functions of $a$. The lowest degree of the Taylor expansion of this functional coefficient at the variable $B$ (respectively, $C$) is one (respectively, two). Therefore, to define a homogeneous limit condition on these coordinates (in which the asymptotics of the coefficients at $B$ and $C$ will be of the same order and the ratio of these coefficients will tend to the prescribed coefficient $\theta/4$ in (\ref{eqfin})), the monomial of degree one in the coefficient at the variable $B$ should vanish. This gives us $p=-\frac{1}{2}$. Then asymptotically $C \approx -2B q$ when $a $ tends to $0$. 

To obtain the desired coefficient $\theta \equiv 4 \frac{C}{B}$ of the limit algebra (\ref{eqfin}), we need to take $q=-\frac{1}{8} \theta$ in (\ref{lamed}). This is the only way to approach each point of the cell $\bar S_2$ from the side of the cell $\bar J^-_{312}$, so we have the incidence coefficient $[\bar J^-_{312}, \bar S_2]=1$. Similar (or sometimes simpler) considerations prove all other incidence coefficients in the formulas of the boundary operators listed in the following subsections.

\subsection{Formulas for boundary operators}

The compact CW-complex $\overline{TD}_2(S^1)$ has only one 0-dimensional cell. Therefore, the boundary operator $\partial_1: C_1 \to C_0$ of the corresponding cellular complex is trivial for any coefficient group.

\begin{proposition}
\label{2to1}
The boundary operator $\partial_2: C_2 \to C_1$ of the mod 2 cellular complex of the CW-complex $\overline{TD}_2(S^1)$ is defined by the following formulas:
\begin{align}
\partial (Y_1) & = \nabla + \bar Y_1 + \bar Y_2 &
\partial (Y_2) & = \nabla + \bar Y_2 + \bar Y_1 \nonumber \\
\partial (U_1) & = \nabla + \bar U_1 + \bar U_2 &
\partial (U_2) & = \nabla + \bar U_2 + \bar U_1 \nonumber \\
\partial (Z) & = \nabla &
\partial (\Theta) & = 0 \nonumber \\
\partial(\bar L_{123}) & = \bar U_1 + \bar Z + \bar Y_2 &
\partial(\bar L_{231}) & = \bar Y_1 + \bar U_2 + \bar Z \nonumber \\
\partial(\bar L_{312}) & = \bar Z + \bar Y_2 + \bar U_2 &
\partial(\bar L_{132}) & = \bar Y_1 + \bar Z + \bar U_2 \nonumber \\
\partial(\bar L_{213}) & = \bar U_1 + \bar Y_2 + \bar Z &
\partial(\bar L_{321}) & = \bar Z + \bar U_2 + \bar Y_2 \nonumber \\
\partial(\bar M_1) & = \bar U_1 + \bar Y_1 + \bar U_2 &
\partial(\bar M_2) & = \bar U_1 + \bar U_2 + \bar Y_1 \nonumber \\
\partial(\bar M_3) & = \bar Y_2 &
\partial (\bar N_1) & = \bar Y_1 + \bar U_2 + \bar Y_2 \nonumber \\
\partial (\bar N_2) & = \bar Y_1 + \bar Y_2 + \bar U_2 &
\partial (\bar N_3) & = \bar U_1 \nonumber \\
\partial(\bar P_1) & = 0 &
\partial (\bar P_2) & = 0 \nonumber \\
\partial (\bar S_1) & = 0 &
\partial(\bar S_2) & = 0 \nonumber \\
\partial(\bar X^+) & = \bar Y_1 + \bar Y_2 &
\partial (\bar X^-) & = \bar Y_1 + \bar Y_2 \nonumber \\
\partial(\bar V_1^+) & = \bar Y_1 + \bar U_1 &
\partial(\bar V_2^+) & = \bar Y_2 + \bar U_2 \nonumber \\
\partial(\bar V_1^-) & = \bar Y_1 + \bar U_1 &
\partial(\bar V_2^-) & = \bar Y_2 + \bar U_2 \nonumber \\
\partial(\bar \Omega) & = 0 \ .\nonumber 
\end{align}
The matrix of this operator is given in page \pageref{page2to1}. \hfill $\Box$
\end{proposition}

\begin{corollary}
\label{corhom1}
The group $H_1(\overline{TD}_2, \Z_2)$ is isomorphic to $\Z_2$ and is generated by the class of the cell $\bar \Theta$.
\end{corollary}

\begin{corollary}
The kernel of the operator $\partial_2 : C_2 \to C_1$ is generated by 23 cycles $($\ref{ker01}$)$--$($\ref{ker23}$):$
\begin{eqnarray}
& &
Y_1+Y_2 \label{ker01} \\ & & 
Y_1+\bar X^+ +Z \label{ker02} \\ & & 
U_1+U_2 \label{ker03} \\ & & 
U_1+Z + \bar M_1 +\bar N_1 + \bar V_2^+ \label{ker04} \\ & & 
\Theta \label{ker05} \\ & & 
\bar L_{123} + \bar L_{213} \label{ker06} \\ & & 
\bar L_{231} + \bar L_{132} \label{ker07} \\ & & 
\bar L_{321}+\bar L_{312} \label{ker08} \\ & & 
\bar L_{231}+\bar L_{312}+\bar X^+ \label{ker09} \\ & & 
\bar L_{123}+\bar L_{312}+\bar M_1 + \bar N_1 + \bar V_2^+ \label{ker10} \\ & & 
\bar M_1+\bar M_2 \label{ker11} \\ & & 
\bar N_1+\bar N_2 \label{ker12} \\ & & 
\bar M_1+\bar N_1 + \bar M_3 + \bar N_3 \label{ker13} \\ & & 
\bar M_3 + \bar N_1 + \bar V_2^+ + \bar X^+ \label{ker14} \\ & & 
\bar M_3+\bar N_3 + \bar X^+ + \bar V_1^+ \label{ker15} \\ & & 
\bar P_1 \label{ker16} \\ & & 
\bar P_2 \label{ker17} \\ & & 
\bar S_1 \label{ker18} \\ & & 
\bar S_2 \label{ker19} \\ & & 
\bar X^+ + \bar X^- \label{ker20} \\ & & 
\bar V_1^+ + \bar V_1^- \label{ker21} \\ & & 
\bar V_2^+ + \bar V_2^- \label{ker22} \\ 
& & \bar \Omega \label{ker23} \qquad \qquad \qquad \qquad \qquad \qquad \qquad \qquad \Box
\end{eqnarray}
\end{corollary}

\begin{proposition}
\label{3to2}
The boundary operator $\partial_3: C_3 \to C_2$ of the mod 2 cell complex $\overline{TD}_2(S^1)$ is defined by following formulas $($both in terms of standard generators of the group $C_2$ and of its subgroup $\ker \partial_1) :$ 
\begin{eqnarray*}
\partial(L_{123}) & = & U_1 + Z + \bar L_{123} + \bar L_{312} = (\ref{ker04}) +(\ref{ker10}) \\
\partial(L_{231}) & = & Y_1 + U_2 + \bar L_{231} + \bar L_{123} = (\ref{ker02}) +(\ref{ker03}) +(\ref{ker04}) +(\ref{ker09}) +(\ref{ker10}) \\
\partial(L_{312}) & = & Z + Y_2 + \bar L_{312} + \bar L_{231} = (\ref{ker01}) +(\ref{ker02}) +(\ref{ker09}) \\
\partial(L_{132}) & = & Y_1 + Z + \bar L_{132} + \bar L_{321} = (\ref{ker02}) +(\ref{ker07}) +(\ref{ker08}) + (\ref{ker09}) \\
\partial(L_{213}) & = & U_1 + Y_2 + \bar L_{213} + \bar L_{132} = (\ref{ker01}) +(\ref{ker02}) +(\ref{ker04}) +(\ref{ker06}) +(\ref{ker07}) +(\ref{ker09}) +(\ref{ker10}) \\
\partial(L_{321}) & = & Z + U_2 + \bar L_{321} + \bar L_{213} = (\ref{ker03}) +(\ref{ker04}) +(\ref{ker06}) + (\ref{ker08}) + (\ref{ker10}) 
\end{eqnarray*}
\begin{eqnarray*}
\partial(M_1) & = & U_1 + Y_1 + \bar M_1 + \bar M_3 = (\ref{ker02}) +(\ref{ker04}) +(\ref{ker14}) \\
\partial(M_2) & = & U_1 + U_2 + \bar M_2 + \bar M_1 = (\ref{ker03}) +(\ref{ker11}) \\
\partial(M_3) & = & Y_2 + U_2 + \bar M_3 + \bar M_2 = (\ref{ker01}) +(\ref{ker02}) +(\ref{ker03}) +(\ref{ker04}) +(\ref{ker11}) +(\ref{ker14}) \\
\partial ( N_1) & = & Y_1 + U_2 + \bar N_1 + \bar N_3 = (\ref{ker02}) +(\ref{ker03}) +(\ref{ker04}) +(\ref{ker13}) +(\ref{ker14}) \\
\partial (N_2) & = & Y_1 + Y_2 + \bar N_2 + \bar N_1 = (\ref{ker01}) +(\ref{ker12}) \\
\partial (N_3) & = & U_1 + Y_2 + \bar N_3 + \bar N_2 = (\ref{ker01}) +(\ref{ker02}) +(\ref{ker04}) +(\ref{ker12}) +(\ref{ker13}) +(\ref{ker14}) \\
\partial(P_1) & = & \Theta + \bar P_1 + \bar P_2 = (\ref{ker05}) +(\ref{ker16}) +(\ref{ker17}) \\
\partial (P_2) & = & \Theta + \bar P_2 + \bar P_1 = (\ref{ker05}) +(\ref{ker16}) +(\ref{ker17}) \\
\partial (S_1) & = & \Theta + \bar S_1 + \bar S_2 = (\ref{ker05}) +(\ref{ker18}) +(\ref{ker19}) \\
\partial(S_2) & = & \Theta + \bar S_2 + \bar S_1 = (\ref{ker05}) +(\ref{ker18}) +(\ref{ker19}) \\
\partial(X^+) & = & Y_1 + Y_2 + \Theta = (\ref{ker01}) +(\ref{ker05}) \\
\partial (X^-) & = & Y_1 + Y_2 = (\ref{ker01}) \\
\partial(V_1^+) & = & Y_1 + U_1 + \bar V_1^+ + \bar V_2^+ = (\ref{ker02}) +(\ref{ker04}) +(\ref{ker13}) +(\ref{ker15}) \\
\partial(V_2^+) & = & Y_2 + U_2 + \bar V_2^+ + \bar V_1^+ = (\ref{ker01}) +(\ref{ker02}) +(\ref{ker03}) +(\ref{ker04}) +(\ref{ker13}) +(\ref{ker15}) \\
\partial(V_1^-) & = & Y_1 + U_1 + \Theta + \bar V_1^- + \bar V_2^- = \\
& & = (\ref{ker02}) +(\ref{ker04}) +(\ref{ker05}) +(\ref{ker13}) +(\ref{ker15}) +(\ref{ker21}) +(\ref{ker22}) \\
\partial(V_2^-) & = & Y_2 + U_2 + \Theta + \bar V_2^- + \bar V_1^- = \\
& & = (\ref{ker01}) +(\ref{ker02}) +(\ref{ker03}) +(\ref{ker04}) +(\ref{ker05}) +(\ref{ker13}) +(\ref{ker15}) +(\ref{ker21}) +(\ref{ker22}) \\
\partial(\Omega) & = & 0 \\
\partial (\bar F_1 ) & = & \bar M_1 + \bar L_{123} + \bar L_{132} + \bar M_3 = (\ref{ker07}) + (\ref{ker09}) + (\ref{ker10}) + (\ref{ker14}) \\
\partial (\bar F_2) & = & \bar M_1 + \bar M_2 + \bar L_{231} + \bar L_{132} = (\ref{ker07}) + (\ref{ker11}) \\
\partial (\bar F_3) & = & \bar L_{213} + \bar M_2 + \bar M_3 + \bar L_{231} = (\ref{ker06}) + (\ref{ker09}) + (\ref{ker10}) + (\ref{ker11}) + (\ref{ker14}) \\
\partial (\bar F_4) & = & \bar L_{312} + \bar L_{321} = (\ref{ker08}) \\
\partial (\bar G_{1234}) & = & \bar L_{123} + \bar N_3 + \bar L_{312} + \bar N_1 + \bar X^- = (\ref{ker10}) + (\ref{ker13}) + (\ref{ker14}) + (\ref{ker20}) \\
\partial (\bar G_{1243}) & = & \bar L_{132} + \bar N_2 + \bar L_{312} + \bar M_3 + \bar V_2^+ = (\ref{ker07}) + (\ref{ker09}) + (\ref{ker12}) + (\ref{ker14}) \\
\partial (\bar G_{1324}) & = & \bar M_1 + \bar V_1^- + \bar S_1 + \bar N_2 + \bar X^- + \bar \Omega + \bar S_2 = \\
& = & (\ref{ker12}) + (\ref{ker13}) + (\ref{ker15}) + (\ref{ker18}) + (\ref{ker19}) + (\ref{ker20}) + (\ref{ker21}) + (\ref{ker23}) \\
\partial (\bar G_{1342}) & = & \bar N_2 + \bar X^- + \bar M_3 + \bar V_2^- + \bar \Omega = (\ref{ker12}) + (\ref{ker14}) + (\ref{ker20}) + (\ref{ker22}) + (\ref{ker23}) \\
\partial (\bar G_{1423}) & = & \bar M_1 + \bar L_{213} + \bar N_3 + \bar V_1^+ + \bar L_{312} = (\ref{ker06}) + (\ref{ker10}) + (\ref{ker14}) + (\ref{ker15}) \\
\partial (\bar G_{1432}) & = & \bar N_3 + \bar L_{213} + \bar N_2 + \bar X^- + \bar L_{321} = \\
& & = (\ref{ker06}) + (\ref{ker08}) + (\ref{ker10}) + (\ref{ker12}) + (\ref{ker13}) + (\ref{ker14}) + (\ref{ker20}) \\
\partial (\bar G_{2134}) & = & \bar L_{123} + \bar M_2 + \bar L_{321} + \bar N_3 + \bar V_1^+ = (\ref{ker08}) + (\ref{ker10}) + (\ref{ker11}) + (\ref{ker14}) + (\ref{ker15}) \\
\partial (\bar G_{2143}) & = & \bar L_{132} + \bar N_1 + \bar L_{321} + \bar N_2 + \bar X^- = (\ref{ker07}) + (\ref{ker08}) + (\ref{ker09}) + (\ref{ker12}) + (\ref{ker20}) \\
\partial (\bar G_{2341}) & = & \bar N_2 + \bar L_{231} + \bar N_1 + X^- + \bar L_{312} = (\ref{ker09}) + (\ref{ker12}) + (\ref{ker20}) \\
\partial (\bar G_{2431}) & = & \bar M_2 + \bar N_1 + \bar X^- + \bar S_1 + \bar \Omega + \bar V_1^- + \bar S_2 = \\
& = & (\ref{ker11}) + (\ref{ker13}) + (\ref{ker15}) + (\ref{ker18}) + (\ref{ker19}) + (\ref{ker20}) + (\ref{ker21}) + (\ref{ker23}) \\
\partial (\bar G_{3142}) & = & \bar M_3 + \bar V_2^- + \bar N_1 + \bar X^- + \bar \Omega = (\ref{ker14}) + (\ref{ker20}) + (\ref{ker22}) + (\ref{ker23}) \\
\partial (\bar G_{3241} ) & = & \bar N_1 + \bar L_{231} + \bar M_3 + \bar V_2^+ + \bar L_{321} = (\ref{ker08}) + (\ref{ker09}) + (\ref{ker14}) 
\end{eqnarray*}
\begin{eqnarray*}
\partial (\bar H_1) & = & \bar M_1 + \bar N_3 + \bar N_2 + \bar M_3 = (\ref{ker12}) + (\ref{ker13}) \\
\partial (\bar H_2) & = & \bar M_1 + \bar M_2 + \bar N_1 + \bar N_2 = (\ref{ker11}) + (\ref{ker12}) \\
\partial (\bar H_3) & = & \bar N_3 + \bar M_2 + \bar M_3 + \bar N_1 = (\ref{ker13}) + (\ref{ker11}) \\
\partial (\bar H_4) & = & \bar N_2 + \bar N_1 = (\ref{ker12}) \\
\partial (\bar I_1) & = & \bar S_1 + \bar P_1 + \bar S_2 = (\ref{ker16}) + (\ref{ker18}) + (\ref{ker19}) \\
\partial (\bar I_2) & = & \bar S_1 + \bar S_2 + \bar P_1 = (\ref{ker16}) + (\ref{ker18}) + (\ref{ker19}) \\
\partial (\bar I_3) & = & \bar P_2 = (\ref{ker17}) \\
\partial (\bar J_{123}^+) & = & \bar X^- + \bar M_1 + \bar N_3 + \bar V_2^+ = (\ref{ker13}) + (\ref{ker14}) + (\ref{ker20}) \\
\partial (\bar J_{123}^-) & = & \bar S_1 + \bar X^+ + \bar M_1 + \bar N_3 + \bar V_2^- = (\ref{ker13}) + (\ref{ker14}) + (\ref{ker18}) + (\ref{ker22}) \\
\partial (\bar J_{132}^+) & = & \bar V_1^- + \bar X^+ + \bar M_1 + \bar N_2 + \bar S_2 = (\ref{ker12}) + (\ref{ker13}) + (\ref{ker15}) + (\ref{ker19}) + (\ref{ker21}) \\
\partial (\bar J_{132}^-) & = & \bar V_1^+ + \bar X^- + \bar M_1 + \bar N_2 = (\ref{ker12}) + (\ref{ker13}) + (\ref{ker15}) + (\ref{ker20}) \\
\partial (\bar J_{213}^+) & = & \bar S_1 + \bar V_2^- + \bar M_2 + \bar N_3 + \bar X^+ = (\ref{ker11}) + (\ref{ker13}) + (\ref{ker14}) + (\ref{ker18}) + (\ref{ker22}) \\
\partial (\bar J_{213}^-) & = & \bar V_2^+ + \bar M_2 + \bar N_3 + \bar X^- = (\ref{ker11}) + (\ref{ker13}) + (\ref{ker14}) + (\ref{ker20}) \\
\partial (\bar J_{231}^+) & = & \bar V_1^+ + \bar M_2 + \bar N_1 + \bar X^- = (\ref{ker11}) + (\ref{ker13}) + (\ref{ker15}) + (\ref{ker20}) \\
\partial (\bar J_{231}^-) & = & \bar V_1^- + \bar S_2 + \bar M_2 + \bar N_1 + \bar X^+ = (\ref{ker11}) + (\ref{ker13}) + (\ref{ker15}) + (\ref{ker19}) + (\ref{ker21}) \\
\partial (\bar J_{312}^+) & = & \bar X^- + \bar V_2^+ + \bar M_3 + \bar N_2 = (\ref{ker12}) + (\ref{ker14}) + (\ref{ker20}) \\
\partial (\bar J_{312}^-) & = & \bar X^+ + \bar V_2^- + \bar M_3 + \bar N_2 + \bar S_2 = (\ref{ker12}) + (\ref{ker14}) + (\ref{ker19}) + (\ref{ker22}) \\
\partial (\bar J_{321}^+) & = & \bar X^+ + \bar M_3 + \bar N_1 + \bar S_2 + \bar V_2^- = (\ref{ker14}) + (\ref{ker19}) + (\ref{ker22}) \\
\partial (\bar J_{321}^-) & = & \bar X^- + \bar M_3 + \bar N_1 + \bar V_2^+ = (\ref{ker14}) + (\ref{ker20}) \\
\partial(\bar K_1^+) & = & \bar V_1^+ + \bar N_3 + \bar N_2 + \bar V_2^+ = (\ref{ker12}) + (\ref{ker14}) + (\ref{ker15}) \\
\partial(\bar K_1^-) & = & \bar V_1^- + \bar P_2 + \bar N_3 + \bar N_2 + \bar V_2^- = (\ref{ker12}) + (\ref{ker14}) + (\ref{ker15}) + (\ref{ker17}) + (\ref{ker21}) + (\ref{ker22}) \\
\partial( \bar K_2^+) & = & \bar V_1^+ + \bar V_2^+ + \bar N_1 + \bar N_3 = (\ref{ker14}) + (\ref{ker15}) \\
\partial( \bar K_2^-) & = & \bar V_1^- + \bar V_2^- + \bar N_1 + \bar N_3 + \bar P_2 = (\ref{ker14}) + (\ref{ker15}) + (\ref{ker17}) + (\ref{ker21}) + (\ref{ker22}) \\
\partial(\bar K_3^+) & = & \bar N_1 + \bar N_2 = (\ref{ker12}) \\
\partial(\bar K_3^-) & = & \bar P_1 + \bar N_1 + \bar N_2 = (\ref{ker12}) + (\ref{ker16}) \\
\partial( \bar W_1^+) & = & \bar P_1 + \bar V_1^+ + \bar V_1^- + \bar \Omega = (\ref{ker16}) + (\ref{ker21}) + (\ref{ker23}) \\
\partial(\bar W_1^-) & = & \bar P_1 + \bar V_1^+ + \bar V_1^- + \bar \Omega = (\ref{ker16}) + (\ref{ker21}) + (\ref{ker23}) \\
\partial(\bar W_2^+) & = & \bar P_2 + \bar V_2^+ + \bar V_2^- + \bar \Omega = (\ref{ker17}) + (\ref{ker22}) + (\ref{ker23}) \\
\partial(\bar W_2^-) & = & \bar P_2 + \bar V_2^+ + \bar V_2^- + \bar \Omega = (\ref{ker17}) + (\ref{ker22}) + (\ref{ker23})
\end{eqnarray*}
The matrix of this operator is given on pages \pageref{3to2lt}--\pageref{3to2rb}. \hfill $\Box$
\end{proposition}

\begin{corollary}
\label{corhom2}
The group $H_2(\overline{TD}_2, Z_2)$ is isomorphic to $\Z_2^2$, its three non-trivial elements are the classes of chains 
\begin{itemize}
\item[a) ] $\bar \Omega \simeq \bar V_1^+ + \bar V_1^- \simeq \bar V_2^+ + \bar V_2^- $,
\item[b) ] $\bar X^+ + \bar X^- $,
\item[c) ] $\bar S_1 \simeq \bar S_2.$ \hfill $\Box$
\end{itemize}
\end{corollary}

\begin{corollary}
The kernel of the operator $\partial_3 : C_3 \to C_2$ is generated by 47 linearly independent chains
\begin{eqnarray}
& & L_{123} + L_{231} + L_{312} + V_1^+ + V_2^+
 \label{kr26} \\
& & L_{132} + L_{312} + \bar G_{2143} + \bar G_{2341} + X^- \label{kr27} \\
& & L_{213} + L_{231} + V_1^+ + X^- + \bar G_{1423} + \bar G_{1234} + \bar J_{123}^+ + \bar G_{2143} + \bar G_{2341} + N_1 + F_4 \label{kr28} \\
& & L_{231} + X^- + V_2^+ + \bar G_{2134} + \bar J_{213}^- + \bar G_{3241} + \bar J_{321}^- \label{kr29} \\
& & L_{312} + L_{321} + V_2^+ + \bar G_{1432} + \bar G_{2341} + \bar K_2^+ \label{kr30} \\
& & M_3 + M_1 + M_2 + X^- \label{kr31} \\
& & M_1 + V_1^+ + \bar K_2^+ + \bar J_{123}^+ + \bar J_{321}^- \label{kr32} \\
& & M_2 + N_1 + V_1^+ + \bar K_2^+ + \bar J_{123}^+ + \bar J_{213}^- \label{kr33} \\
& & N_1 + X^- + V_2^+ + \bar K_2^+ \label{kr34} \\
& & N_2 + \bar K_3^+ + X^- \label{kr35} \\
& & N_3 + V_1^+ + \bar K_1^+ + X^- \label{kr36} \\
& & P_1 + P_2 \label{kr37} \\
& & P_2 + S_2 + \bar I_1 + \bar I_3 \label{kr38} \\
& & S_1 + S_2 \label{kr39} \\
& & S_2 + \bar I_1 + \bar I_3 + V_1^+ + V_1^- + \bar W_1^+ + \bar W_2^+ \label{kr40} \\
& & X^+ + X^- + \bar K_3^+ + \bar K_3^- + \bar I_3 + V_1^+ + V_1^- + \bar W_1^+ + \bar W_2^+ \label{kr41} \\
& & V_1^+ + V_2^+ + V_1^- + V_2^- \label{kr42} \\
& & \Omega \label{kr43} \\
& & \bar F_1 + \bar G_{1234} + \bar G_{1243} + \bar J_{123}^+ + \bar K_3^+ \label{kr44} \\
& & \bar F_2 + F_4 + \bar J_{213}^- + \bar J_{123}^+ + \bar G_{2143} + \bar G_{2341} \label{kr45} \\
& & \bar F_3 + \bar G_{1423} + \bar K_1^+ + \bar G_{2341} + \bar J_{321}^- + \bar J_{213}^- + \bar J_{123}^+ \label{kr46} \\
& & \bar F_4 + \bar G_{3241} + \bar G_{2341} + \bar K_3^+ + \bar J_{321}^- \label{kr47} \\
& & \bar G_{1234} + \bar J_{213}^- + \bar G_{2134} + \bar G_{3241} + \bar G_{2341} + \bar K_1^+ + \bar J_{321}^- \label{kr48} \\
& & \bar G_{1243} + \bar G_{2143} + \bar G_{2341} + \bar G_{3241} + \bar K_3^+ \label{kr49} \\
& & \bar G_{1324} + \bar G_{2431} + \bar H_2 \label{kr50} \\
& & \bar G_{1342} + \bar G_{3142} + \bar K_3^+ \label{kr51} \\
& & \bar G_{1423} + \bar K_2^+ + \bar G_{1432} + \bar J_{123}^+ + \bar G_{3241} + \bar G_{2341} + \bar J_{321}^- \label{kr52} \\
& & \bar G_{2431} + \bar H_3 + \bar K_3^+ + \bar J_{132}^- + \bar I_1 + \bar W_1^+ + \bar J_{123}^+ + \bar J_{321}^- \label{kr53} \\
& & \bar G_{3142} + \bar W_2^+ + \bar I_3 + \bar J_{321}^- \label{kr54} \\
& & \bar H_3 + \bar H_1 + \bar H_2 \label{kr55} \\
& & \bar H_1 + \bar J_{132}^- + \bar J_{321}^- + \bar K_2^+ \label{kr56} \\
& & \bar H_2 + \bar J_{213}^- + \bar J_{123}^+ + \bar K_3^+ \label{kr57} \\
& & \bar H_4 + \bar K_3^+ \label{kr58} 
\end{eqnarray}
\begin{eqnarray}
& & \bar I_1 + \bar I_2 \label{kr59} \\
& & \bar I_2 + \bar J_{123}^+ + \bar J_{123}^- + \bar J_{132}^+ + \bar J_{132}^- + 
\bar I_3 + \bar W_1^+ + \bar W_2^+ \label{kr60} \\
& & \bar I_3 + \bar J_{132}^+ + \bar J_{132}^- + \bar J_{321}^+ + \bar J_{321}^- + \bar K_2^+ + \bar K_2^- \label{kr61} \\
& & \bar J_{123}^- + \bar J_{123}^+ + \bar J_{213}^+ + \bar J_{213}^- \label{kr62} \\
& & \bar J_{123}^+ + \bar K_1^+ + \bar J_{132}^- \label{kr63} \\
& & \bar J_{132}^+ + \bar J_{132}^- + \bar J_{231}^+ + \bar J_{231}^- \label{kr64} \\
& & \bar J_{213}^- + \bar J_{231}^+ + \bar K_2^+ \label{kr65} \\
& & \bar J_{312}^+ + \bar J_{321}^- + \bar K_3^+ \label{kr66} \\
& & \bar J_{312}^- + \bar J_{321}^+ + \bar K_3^+ \label{kr67} \\
& & \bar K_1^- + \bar K_2^- + \bar K_3^+ \label{kr68} \\
& & \bar K_1^+ + \bar K_2^+ + \bar K_3^+ \label{kr69} \\
& & \bar K_2^+ + \bar K_2^- + \bar K_3^+ + \bar K_3^- + \bar W_1^+ + \bar W_2^+ \label{kr70} \\
& & \bar W_1^+ + \bar W_1^- \label{kr71} \\
& & \bar W_2^+ + \bar W_2^- \label{kr72} \qquad \qquad \qquad \qquad \qquad \qquad \qquad \Box
\end{eqnarray}
\end{corollary}

\begin{proposition}
\label{4to3}
The boundary operator $\partial_4: C_4 \to C_3$ of the mod 2 cell complex $\overline{TD}_2(S^1)$ is defined by following formulas $($both in the terms of three-dimensional cells and the basic three-dimensional cycles $($\ref{kr26}$)$--$($\ref{kr72}$))$:
\begin{eqnarray*}
\partial (F_1 ) & = & M_1 + L_{123} + L_{132} + \bar F_1 + \bar F_4 =\\
& & (\ref{kr26})+ (\ref{kr27})+ (\ref{kr29})+ (\ref{kr32})+ (\ref{kr44})+ (\ref{kr47})+ (\ref{kr48})+ (\ref{kr49}) + (\ref{kr69})\nonumber \\
\partial(F_2) & = &M_1 + M_2 + L_{231} + \bar F_1 + \bar F_2 = \\
& & (\ref{kr29})+ (\ref{kr32})+ (\ref{kr33})+ (\ref{kr34})+ (\ref{kr44})+ (\ref{kr45})+ (\ref{kr47})+ (\ref{kr48})+ (\ref{kr49})+ (\ref{kr69}) \nonumber \\
\partial(F_3) & = &L_{213} + M_2 + M_3 + \bar F_3 + \bar F_2 =\\
& & = (\ref{kr28})+ (\ref{kr29})+ (\ref{kr31})+ (\ref{kr32})+ (\ref{kr34})+ (\ref{kr45})+ (\ref{kr46})+ (\ref{kr48}) \nonumber \\
\partial (F_4) & = &L_{312} + L_{321} + M_3 + \bar F_4 + \bar F_3 = \\
& & = (\ref{kr30})+ (\ref{kr31})+ (\ref{kr32})+ (\ref{kr33})+ (\ref{kr34})+ (\ref{kr46})+ (\ref{kr47})+ (\ref{kr52}) + (\ref{kr69}) \nonumber \\
\partial(G_{1234}) & = &L_{123} + N_3 + L_{312} + \bar G_{1234} + \bar G_{2341} = 
 (\ref{kr26})+ (\ref{kr29})+ (\ref{kr36})+ (\ref{kr48}) \\
\partial (G_{1243}) & = & L_{132} + N_2 + L_{312} + \bar G_{1243} + \bar G_{3241} = (\ref{kr27})+ (\ref{kr35})+ (\ref{kr49}) \\
\partial (G_{1324}) & = & M_1 + V_1^- + S_1 + \bar G_{1324} + \bar G_{1342} = \\
& & = (\ref{kr32})+ (\ref{kr39})+ (\ref{kr40})+ (\ref{kr50})+ (\ref{kr51})+ (\ref{kr53})+ (\ref{kr54})+ (\ref{kr55})+ (\ref{kr56}) \nonumber \\
\partial (G_{1342}) & = & N_2 + X^- + \bar G_{1342} + \bar G_{3142} = (\ref{kr35})+ (\ref{kr51}) \\
\partial (G_{1423}) & = &M_1 + L_{213} + N_3 + V_1^+ + \bar G_{1423} + \bar G_{1243} = \\
& & = (\ref{kr28})+ (\ref{kr29})+ (\ref{kr32})+ (\ref{kr34})+ (\ref{kr36})+ (\ref{kr47})+ (\ref{kr48})+ (\ref{kr49}) \nonumber \\
\partial (G_{1432}) & = &N_3 + L_{213} + N_2 + X^- + \bar G_{1432} + \bar G_{2143} = \\
& & = (\ref{kr28})+ (\ref{kr29})+ (\ref{kr34})+ (\ref{kr35})+ (\ref{kr36})+ (\ref{kr47})+ (\ref{kr48})+ (\ref{kr52}) \\
\partial (G_{2134}) & = &L_{123} + M_2 + L_{321} + \bar G_{2134} + \bar G_{1423} = \\
& & = (\ref{kr26})+ (\ref{kr29})+ (\ref{kr30})+ (\ref{kr33})+ (\ref{kr34})+ (\ref{kr52}) 
\end{eqnarray*}
\begin{eqnarray*}
\partial (G_{2143}) & = &L_{132} + N_1 + L_{321} + \bar G_{2143} + \bar G_{1432} = (\ref{kr27})+ (\ref{kr30})+ (\ref{kr34}) \\
\partial (G_{2341}) & = &N_2 + L_{231} + N_1 + X^- + \bar G_{2341} + \bar G_{1234} = \\
& & = (\ref{kr29})+ (\ref{kr34})+ (\ref{kr35})+ (\ref{kr48})+ (\ref{kr69}) \\
\partial (G_{2431}) & = &N_3 + M_2 + N_1 + X^- + S_1 + S_2 + \Omega + \bar G_{2431} + \bar G_{1324} = \\
& & = (\ref{kr33})+ (\ref{kr36})+ (\ref{kr39})+ (\ref{kr43})+ (\ref{kr50})+ (\ref{kr57})+ (\ref{kr69}) \\
\partial (G_{3142}) & = & M_3 + V_2^- + S_2 + \bar G_{3142} + \bar G_{2431} = \\
& & = (\ref{kr31})+ (\ref{kr32})+ (\ref{kr33})+ (\ref{kr34})+ (\ref{kr40})+ (\ref{kr42})+ (\ref{kr53})+ (\ref{kr54})+ \\
& & + (\ref{kr55})+ (\ref{kr56})+ (\ref{kr57}) \\
\partial (G_{3241} ) & = & N_1 + L_{231} + M_3 + V_2^+ + \bar G_{3241} + \bar G_{2134} = \\
& & = (\ref{kr29})+ (\ref{kr31})+ (\ref{kr32})+ (\ref{kr33}) \\
\partial (H_1) & = &M_1 + N_3 + N_2 + \bar H_1 + \bar H_4 = \\
& & = (\ref{kr32})+ (\ref{kr35})+ (\ref{kr36})+ (\ref{kr56})+ (\ref{kr58})+ (\ref{kr63})\\
\partial (H_2) & = &M_1 + M_2 + N_1 + \bar H_2 + \bar H_1 = \\
& & = (\ref{kr32})+ (\ref{kr33})+ (\ref{kr56})+ (\ref{kr57})+ (\ref{kr63})+ (\ref{kr69})\\
\partial (H_3) & = &N_3 + M_2 + M_3 + \bar H_3 + \bar H_2 = \\
& & = (\ref{kr31})+ (\ref{kr32})+ (\ref{kr36})+ (\ref{kr55})+ (\ref{kr56})+ (\ref{kr63}) \\
\partial (H_4) & = & N_2 + N_1 + M_3 + \bar H_4 + \bar H_3 = \\
& & = (\ref{kr31})+ (\ref{kr32})+ (\ref{kr33})+ (\ref{kr35})+ (\ref{kr55})+ (\ref{kr56})+ (\ref{kr57})+ (\ref{kr58})+ (\ref{kr63})+ (\ref{kr69})\\
\partial (I_1) & = & S_1 + P_1 + \bar I_1 + \bar I_3 = (\ref{kr37})+ (\ref{kr38})+ (\ref{kr39})\\
\partial (I_2) & = & S_1 + S_2 + \bar I_2 + \bar I_1 = (\ref{kr39})+ (\ref{kr59}) \\
\partial (I_3) & = & P_2 + S_2 + \bar I_3 + \bar I_2 = (\ref{kr38})+ (\ref{kr59}) \\
\partial (J_{123}^+) & = & X^- + M_1 + N_3 + \bar J_{123}^+ + \bar J_{312}^+ = (\ref{kr32})+ (\ref{kr36})+ (\ref{kr66})+ (\ref{kr69}) \\
\partial (J_{123}^-) & = & S_1 + X^+ + M_1 + N_3 + \bar J_{123}^- + \bar J_{312}^- =\\
& & = (\ref{kr32})+ (\ref{kr36})+ (\ref{kr39})+ (\ref{kr40})+ (\ref{kr41})+ (\ref{kr59})+ (\ref{kr60})+ (\ref{kr61})+ \\
& & + (\ref{kr67})+ (\ref{kr69})+ (\ref{kr70}) \\
\partial (J_{132}^+) & = & V_1^- + X^+ + M_1 + N_2 + \bar J_{132}^+ + \bar J_{321}^+ = \\
& & = (\ref{kr32})+ (\ref{kr35})+ (\ref{kr41})+ (\ref{kr61})+ (\ref{kr63})+ (\ref{kr69})+ (\ref{kr70}) \\
\partial (J_{132}^-) & = & V_1^+ + X^- + M_1 + N_2 + \bar J_{132}^- + \bar J_{321}^- = \\
& & = (\ref{kr32})+ (\ref{kr35})+ (\ref{kr63})+ (\ref{kr69}) \\
\partial (J_{213}^+) & = & S_1 + V_2^- + M_2 + N_3 + \bar J_{213}^+ + \bar J_{132}^+ = \\
& & = (\ref{kr33})+ (\ref{kr34})+ (\ref{kr36})+ (\ref{kr39})+ (\ref{kr40})+ (\ref{kr42})+ (\ref{kr59})+ (\ref{kr60})+ (\ref{kr62})+ (\ref{kr63}) \\
\partial (J_{213}^-) & = & V_2^+ + M_2 + N_3 + \bar J_{213}^- + \bar J_{132}^- = (\ref{kr33})+ (\ref{kr34})+ (\ref{kr36})+ (\ref{kr63}) \\
\partial (J_{231}^+) & = & V_1^+ + M_2 + N_1 + \bar J_{231}^+ + \bar J_{123}^+ = (\ref{kr33})+ (\ref{kr65}) \\
\partial (J_{231}^-) & = & V_1^- + S_2 + M_2 + N_1 + \bar J_{231}^- + \bar J_{123}^- = \\
& & = (\ref{kr33})+ (\ref{kr40})+ (\ref{kr59})+ (\ref{kr60})+ (\ref{kr64})+ (\ref{kr65}) \\
\partial (J_{312}^+) & = & X^- + V_2^+ + M_3 + N_2 + \bar J_{312}^+ + \bar J_{231}^+ = \\
& & = (\ref{kr31})+ (\ref{kr32})+ (\ref{kr33})+ (\ref{kr34})+ (\ref{kr35})+ (\ref{kr65})+ (\ref{kr66}) 
\end{eqnarray*}
\begin{eqnarray*}
\partial (J_{312}^-) & = & X^+ + V_2^- + M_3 + N_2 + \bar J_{312}^- + \bar J_{231}^- = \\
& & = (\ref{kr31})+ (\ref{kr32})+ (\ref{kr33})+ (\ref{kr34})+ (\ref{kr35})+ (\ref{kr41})+ (\ref{kr42})+ (\ref{kr61})+ (\ref{kr64})+\\
& & + (\ref{kr65})+ (\ref{kr67})+ (\ref{kr70}) \\
\partial (J_{321}^+) & = & X^+ + M_3 + N_1 + S_2 + \bar J_{321}^+ + \bar J_{213}^+ = \\
& & = (\ref{kr31})+ (\ref{kr32})+ (\ref{kr33})+ (\ref{kr40})+ (\ref{kr41})+ (\ref{kr59})+ (\ref{kr60})+ (\ref{kr61})+ (\ref{kr62})+ (\ref{kr70}) \\
\partial (J_{321}^-) & = & X^- + M_3 + N_1 + \bar J_{321}^- + \bar J_{213}^- = (\ref{kr31})+ (\ref{kr32})+ (\ref{kr33}) \\
\partial( K_1^+) & = & V_1^+ + N_3 + N_2 + \bar K_1^+ + \bar K_3^+ = (\ref{kr35})+ (\ref{kr36}) \\
\partial( K_1^-) & = & V_1^- + P_2 + N_3 + N_2 + \bar K_1^- + \bar K_3^- = \\
& & = (\ref{kr35})+ (\ref{kr36})+ (\ref{kr38})+ (\ref{kr40})+ (\ref{kr68})+ (\ref{kr69})+ (\ref{kr70}) \\
\partial( K_2^+) & = & V_1^+ + V_2^+ + N_1 + N_3 + \bar K_2^+ + \bar K_1^+ = (\ref{kr34})+ (\ref{kr36}) \\
\partial( K_2^-) & = & V_1^- + V_2^- + N_1 + N_3 + \bar K_2^- + \bar K_1^- = \\
& & = (\ref{kr34})+ (\ref{kr36})+ (\ref{kr42})+ (\ref{kr68})+ (\ref{kr69}) \\
\partial( K_3^+) & = & V_2^+ + N_1 + N_2 + \bar K_3^+ + \bar K_2^+ = (\ref{kr34})+ (\ref{kr35}) \\
\partial( K_3^-) & = & P_1 + V_2^- + N_1 + N_2 + \bar K_3^- + \bar K_2^- = \\
& & = (\ref{kr34})+ (\ref{kr35})+ (\ref{kr37})+ (\ref{kr38})+ (\ref{kr40})+ (\ref{kr42})+ (\ref{kr70}) \\
\partial( W_1^+) & = & P_1 + V_1^+ + V_1^- + \bar W_1^+ + \bar W_2^+ = (\ref{kr37})+ (\ref{kr38})+ (\ref{kr40}) \\
\partial( W_1^-) & = & P_1 + V_1^+ + V_1^- + \Omega + \bar W_2^- + \bar W_1^- = \\
& & = (\ref{kr37})+ (\ref{kr38})+ (\ref{kr40})+ (\ref{kr43})+ (\ref{kr71})+ (\ref{kr72}) \\
\partial( W_2^+) & = & P_2 + V_2^+ + V_2^- + \Omega + \bar W_1^+ + \bar W_2^+ = (\ref{kr38})+ (\ref{kr40})+ (\ref{kr42})+ (\ref{kr43}) \\
\partial( W_2^-) & = & P_2 + V_2^+ + V_2^- + \bar W_2^- + \bar W_1^- =\\ 
& & = (\ref{kr38})+ (\ref{kr40})+ (\ref{kr42})+ (\ref{kr71})+ (\ref{kr72}) \\
\partial(\bar B) & = & \bar H_1 + \bar H_2 +\bar H_3 = (\ref{kr55}) \\
 \partial(\bar C_{12}) & = & \bar F_1 + \bar H_1 + \bar G_{1234} + \bar G_{1243} + \bar H_4 + \bar J_{312}^+ = \\
& & = (\ref{kr44})+ (\ref{kr56})+ (\ref{kr58})+ (\ref{kr63})+ (\ref{kr66})+ (\ref{kr69}) \\
 \partial(\bar C_{13}) & = & \bar H_1 + \bar J_{123}^+ + \bar G_{1342} + \bar H_4 + \bar I_3 + \bar W_2^+ + \bar J_{312}^- + \bar J_{321}^+ = \\
& & = (\ref{kr51})+ (\ref{kr54})+ (\ref{kr56})+ (\ref{kr58})+ (\ref{kr63})+ (\ref{kr67})+ (\ref{kr69}) \\
 \partial(\bar C_{14}) & = & \bar H_1 + \bar G_{1324} + \bar J_{123}^- + \bar J_{132}^+ + \bar I_3 + \bar W_2^- + \bar H_4 + \bar J_{321}^- = \\
& & = (\ref{kr50})+ (\ref{kr53})+ (\ref{kr55})+ (\ref{kr58})+ (\ref{kr59})+ (\ref{kr60})+ (\ref{kr72}) \\
 \partial(\bar C_{15}) & = & \bar G_{1423} + \bar G_{1432} + \bar J_{132}^- + \bar F_4 = (\ref{kr47})+ (\ref{kr52})+ (\ref{kr63})+ (\ref{kr69}) \\
 \partial(\bar C_{21}) & = & \bar F_1 + \bar H_2 + \bar G_{2134} +\bar G_{2143} + \bar H_1 + \bar J_{132}^- = \\
& & = (\ref{kr44})+ (\ref{kr48})+ (\ref{kr49})+ (\ref{kr56})+ (\ref{kr57})+ (\ref{kr69}) \\
 \partial(\bar C_{23}) & = & \bar H_1 + \bar F_2 + \bar H_2 + \bar G_{2341} + \bar J_{123}^+ + \bar G_{1243} = \\
& & = (\ref{kr45})+ (\ref{kr47})+ (\ref{kr49})+ (\ref{kr56})+ (\ref{kr57})+ (\ref{kr63})+ (\ref{kr69}) \\
 \partial(\bar C_{24}) & = & \bar H_1 + \bar H_2 + \bar J_{123}^- + \bar J_{132}^+ + \bar J_{231}^+ + \bar I_1 + \bar W_1^+ + \bar G_{1342} = \\
& & = (\ref{kr51})+ (\ref{kr54})+ (\ref{kr56})+ (\ref{kr57})+ (\ref{kr59})+ (\ref{kr60})+ (\ref{kr65}) \\
\partial(\bar C_{25}) & = & \bar G_{2431} + \bar J_{132}^- + \bar J_{231}^- + \bar I_1 + \bar W_1^- + \bar J_{132}^+ = \\
& & = (\ref{kr53})+ (\ref{kr55})+ (\ref{kr56})+ (\ref{kr57})+ (\ref{kr64})+ (\ref{kr65})+ (\ref{kr71}) 
\end{eqnarray*}
\begin{eqnarray*}
 \partial(\bar C_{31}) & = & \bar H_2 + \bar H_3 + \bar G_{3142} + \bar J_{213}^+ + \bar J_{132}^- + \bar J_{231}^- + \bar I_2 + \bar W_1^- = \\
& & = (\ref{kr54})+ (\ref{kr55})+ (\ref{kr56})+ (\ref{kr60})+ (\ref{kr62})+ (\ref{kr64})+ (\ref{kr65})+ (\ref{kr71}) \\
 \partial(\bar C_{32}) & = & \bar H_2 + \bar F_2 + \bar H_3 + \bar G_{3241} + \bar J_{213}^- + \bar G_{2143} = \\
& & = (\ref{kr45})+ (\ref{kr47})+ (\ref{kr55})+ (\ref{kr56})+ (\ref{kr63})+ (\ref{kr69}) \\
\partial(\bar C_{34}) & = & \bar G_{1423} + \bar H_2 + \bar F_3 + \bar H_3 + \bar J_{231}^+ + \bar G_{2341} = \\
& & = (\ref{kr46})+ (\ref{kr55})+ (\ref{kr56})+ (\ref{kr63})+ (\ref{kr65}) \\
 \partial(\bar C_{35}) & = & \bar G_{1324} + \bar J_{231}^- + \bar J_{132}^+ + \bar J_{231}^+ + \bar I_2 + \bar W_1^+ = \\
& & = (\ref{kr50})+ (\ref{kr53})+ (\ref{kr55})+ (\ref{kr56})+ (\ref{kr59})+ (\ref{kr63})+ (\ref{kr64})+ (\ref{kr69}) \\
 \partial(\bar C_{41}) & = & \bar G_{2431} + \bar H_4 + \bar J_{213}^+ + \bar J_{312}^+ + \bar I_3 + \bar H_3 + \bar J_{231}^- + \bar W_2^+ = \\
& & = (\ref{kr53})+ (\ref{kr58})+ (\ref{kr59})+ (\ref{kr60})+ (\ref{kr62})+ (\ref{kr63})+ (\ref{kr64})+ (\ref{kr65})+ (\ref{kr66})+ (\ref{kr69}) \\
 \partial(\bar C_{42}) & = & \bar H_3 + \bar H_4 + \bar J_{213}^- + \bar J_{312}^- + \bar J_{321}^+ + \bar G_{3142} + W_2^- + \bar I_3 = \\
& & = (\ref{kr54})+ (\ref{kr55})+ (\ref{kr56})+ (\ref{kr57})+ (\ref{kr58})+ (\ref{kr63})+ (\ref{kr67})+ (\ref{kr69})+ (\ref{kr72}) \\
 \partial(\bar C_{43}) & = & \bar G_{1432} + \bar H_3 + \bar F_3 + \bar H_4 + \bar J_{321}^- + \bar G_{3241} = \\
& & = (\ref{kr46})+ (\ref{kr52})+ (\ref{kr55})+ (\ref{kr56})+ (\ref{kr57})+ (\ref{kr58})+ (\ref{kr63}) \\
\partial(\bar C_{45}) & = & \bar G_{1234} + \bar G_{2134} + \bar F_4 + \bar J_{231}^+ = (\ref{kr47})+ (\ref{kr48})+ (\ref{kr65})+ (\ref{kr69}) \\
\partial(\bar C_{51}) & = & \bar G_{2341} + \bar G_{3241} + \bar J_{312}^+ + \bar F_4 = (\ref{kr47})+ (\ref{kr66}) \\
\partial (\bar C_{52}) & = & \bar G_{3142} + \bar J_{312}^- + \bar J_{321}^+ + \bar J_{312}^+ + \bar I_3 + \bar W_2^+ = (\ref{kr54})+ (\ref{kr66})+ (\ref{kr67}) \\
 \partial (\bar C_{53}) & = & \bar G_{1342} + \bar J_{321}^- + \bar J_{312}^- + \bar J_{321}^+ + \bar I_3 + \bar W_2^- = (\ref{kr51})+ (\ref{kr54})+ (\ref{kr67})+ (\ref{kr72}) \\
 \partial(\bar C_{54}) & = & \bar G_{1243} + \bar G_{2143} + \bar F_4 + \bar J_{321}^- = (\ref{kr47})+ (\ref{kr49}) \\
 \partial(\bar D_{12}^+) & = & \bar J_{213}^- + \bar K_3^+ + \bar H_1 + \bar H_2 + \bar J_{312}^+ = (\ref{kr56})+ (\ref{kr57})+ (\ref{kr63})+ (\ref{kr66})+ (\ref{kr69}) \\
 \partial(\bar D_{12}^-) & = & \bar I_1 + \bar J_{213}^+ + \bar K_3^- + \bar H_1 + \bar H_2 + \bar J_{312}^- = \\
& & = (\ref{kr56})+ (\ref{kr57})+ (\ref{kr59})+ (\ref{kr60})+ (\ref{kr61})+ (\ref{kr62})+ (\ref{kr63})+ (\ref{kr67})+ (\ref{kr69})+ (\ref{kr70}) \\
 \partial(\bar D_{13}^+) & = & \bar J_{123}^- + \bar J_{213}^+ + \bar J_{312}^- + \bar H_1 + \bar H_3 + \bar J_{321}^+ = (\ref{kr55})+ (\ref{kr57})+ (\ref{kr62})+ (\ref{kr67}) \\
 \partial(\bar D_{13}^-) & = & \bar J_{123}^+ + \bar J_{213}^- + \bar J_{312}^+ + \bar H_1 + \bar H_3 + \bar J_{321}^- = (\ref{kr55})+ (\ref{kr57})+ (\ref{kr66}) \\
 \partial(\bar D_{14}^+) & = & \bar J_{132}^- + \bar K_2^+ + \bar J_{312}^+ + \bar H_1 + \bar H_4 = (\ref{kr56})+ (\ref{kr58})+ (\ref{kr66}) \\
 \partial(\bar D_{14}^-) & = & \bar J_{132}^+ + \bar K_2^- + \bar J_{312}^- + \bar H_1 + \bar H_4 + \bar I_3 = (\ref{kr56})+ (\ref{kr58})+ (\ref{kr61})+ (\ref{kr67}) \\
 \partial(\bar D_{23}^+) & = & \bar J_{123}^+ + \bar J_{321}^- + \bar H_2 + \bar H_3 + \bar K_3^+ = (\ref{kr55})+ (\ref{kr56})+ (\ref{kr63})+ (\ref{kr69}) \\
 \partial(\bar D_{23}^-) & = & \bar J_{123}^- + \bar I_2 + \bar J_{321}^+ + \bar H_2 + \bar H_3 + \bar K_3^- = \\
& & = (\ref{kr55})+ (\ref{kr56})+ (\ref{kr60})+ (\ref{kr61})+ (\ref{kr63})+ (\ref{kr69})+ (\ref{kr70}) \\
 \partial(\bar D_{24}^+) & = & \bar J_{132}^+ + \bar J_{231}^- + \bar J_{321}^+ + \bar H_2 + \bar H_4 + \bar J_{312}^- = \\
& & = (\ref{kr57})+ (\ref{kr58})+ (\ref{kr63})+ (\ref{kr64})+ (\ref{kr65})+ (\ref{kr67})+ (\ref{kr69}) \\
 \partial(\bar D_{24}^-) & = & \bar J_{132}^- + \bar J_{231}^+ + \bar J_{321}^- + \bar H_2 + \bar H_4 + \bar J_{312}^+ = \\
& & = (\ref{kr57})+ (\ref{kr58})+ (\ref{kr63})+ (\ref{kr65})+ (\ref{kr66})+ (\ref{kr69}) \\
 \partial(\bar D_{34}^+) & = & \bar K_1^+ + \bar J_{231}^+ + \bar H_3 + \bar H_4 + J_{321}^- = \\
& & = (\ref{kr55})+ (\ref{kr56})+ (\ref{kr57})+ (\ref{kr58})+ (\ref{kr63})+ (\ref{kr65}) \\
 \partial(\bar D_{34}^-) & = & \bar K_1^- + \bar J_{231}^- + \bar H_3 + \bar H_4 + J_{321}^+ + I_3 = \\
& & = (\ref{kr55})+ (\ref{kr56})+ (\ref{kr57})+ (\ref{kr58})+ (\ref{kr61})+ (\ref{kr63})+ (\ref{kr64})+ (\ref{kr65})+ (\ref{kr68})+ (\ref{kr69}) 
\end{eqnarray*}
\begin{eqnarray*}
 \partial(\bar E_+) & = & \bar K_1^+ + \bar K_2^+ + \bar K_3^+ = (\ref{kr69}) \\
 \partial (\bar E_1) & = & \bar K_1^+ + \bar K_2^- + \bar K_3^- + \bar W_1^+ + \bar W_2^-
= (\ref{kr69})+ (\ref{kr70})+ (\ref{kr72}) \\
 \partial(\bar E_2) & = & \bar K_1^- + \bar K_2^+ + \bar K_3^- + \bar W_1^- + \bar W_2^+ =
 (\ref{kr68})+ (\ref{kr70})+ (\ref{kr71}) \\
 \partial(\bar E_3) & = & \bar K_1^- + \bar K_2^- + \bar K_3^+ + \bar W_2^- + W_2^+ = (\ref{kr68})+ (\ref{kr72}) \qquad \qquad \qquad \qquad \quad \Box
\end{eqnarray*}
\end{proposition}

\begin{corollary}
\label{corhom3}
The group $H_3(\overline{TD}_2(S^1), \Z_2)$ is isomorphic to $\Z_2^2$. Its three non-zero elements are defined by the classes of cycles 
\begin{itemize}
\item[a) ] $ (\ref{kr53}) \simeq (\ref{kr67})$,
\item[b) ] $ (\ref{kr41})\simeq (\ref{kr54}) \simeq (\ref{kr60}) \simeq (\ref{kr61}) \simeq (\ref{kr62}) \simeq (\ref{kr64}) $,
\item[c) ] $ (\ref{kr26}) \simeq (\ref{kr30}) \simeq (\ref{kr31}) \simeq (\ref{kr33}) \simeq (\ref{kr34}) \simeq (\ref{kr35}) \simeq (\ref{kr36}) \simeq (\ref{kr44}) \simeq (\ref{kr45}) \simeq (\ref{kr46}) \simeq (\ref{kr47}) \simeq (\ref{kr50}) \simeq (\ref{kr51}) \simeq (\ref{kr56}) \simeq (\ref{kr57}) \simeq (\ref{kr63}) \simeq (\ref{kr65}) \simeq (\ref{kr66})$.
\end{itemize}
All remaining cycles from the list $ (\ref{kr26})$--$ (\ref{kr72})$ are homologous to zero.
\hfill $\Box$
\end{corollary}

\begin{proposition}
\label{5to4}
The boundary operator $\partial_5: C_5 \to C_4$ of the mod 2 cell complex $\overline{TD}_2(S^1)$ is defined by following formulas:
\label{5 to 4} 
\begin{eqnarray*}
 \partial(B) & = & H_1 + H_2 + H_3 + H_4 \\
 \partial(C_{12}) & = & F_1 + H_1 + G_{1234} + G_{1243} + \bar C_{12} + \bar C_{51} \\
 \partial(C_{13}) & = & H_1 + J_{123}^+ + G_{1342} + \bar C_{13} + \bar C_{52} \\
 \partial(C_{14}) & = & H_1 + G_{1324} + J_{123}^- + J_{132}^+ + \bar C_{14} + \bar C_{53} \\
 \partial(C_{15}) & = & G_{1423} + G_{1432} + J_{132}^- + \bar C_{15} + \bar C_{54}\\
 \partial(C_{21}) & = & F_1 + H_2 + G_{2134} + G_{2143} + \bar C_{21} + \bar C_{15} \\
 \partial(C_{23}) & = & H_1 + F_2 + H_2 + G_{2341} + J_{123}^+ + \bar C_{23} + \bar C_{12} \\
 \partial(C_{24}) & = & H_1 + H_2 + J_{123}^- + J_{132}^+ + J_{231}^+ + I_1 + W_1^+ + \bar C_{24} + \bar C_{13} \\
\partial(C_{25}) & = & H_1 + G_{2431} + J_{132}^- + J_{231}^- + I_1 + W_1^- + \bar C_{25} + \bar C_{14} \\
 \partial(C_{31}) & = & H_3 + G_{3142} + J_{213}^+ + I_2 + \bar C_{31} + \bar C_{25} \\
 \partial(C_{32}) & = & H_2 + F_2 + H_3 + G_{3241} + J_{213}^- + \bar C_{32} + \bar C_{21} \\
\partial(C_{34}) & = & G_{1423} + H_2 + F_3 + H_3 + J_{231}^+ + \bar C_{34} + \bar C_{23} \\
 \partial(C_{35}) & = & G_{1324} + H_2 + I_2 + J_{231}^- + \bar C_{35} + \bar C_{24}\\
 \partial(C_{41}) & = & G_{2431} + H_4 + J_{213}^+ + J_{312}^+ + I_3 + W_2^+ + \bar C_{41} + \bar C_{35} \\
 \partial(C_{42}) & = & H_3 + H_4 + J_{213}^- + J_{312}^- + J_{321}^+ + I_3 + W_2^- + \bar C_{42} + \bar C_{31} \\
 \partial(C_{43}) & = & G_{1432} + H_3 + F_3 + H_4 + J_{321}^- + \bar C_{43} + \bar C_{32} \\
\partial(C_{45}) & = & G_{1234} + G_{2134} + H_3 + F_4 + \bar C_{45} + \bar C_{34} \\
\partial(C_{51}) & = & G_{2341} + G_{3241} + J_{312}^+ + \bar C_{51} + \bar C_{45} \\
\partial (C_{52}) & = & H_4 + G_{3142} + J_{312}^- + J_{321}^+ + \bar C_{52} + \bar C_{41} \\
 \partial (C_{53}) & = & G_{1342} + H_4 + J_{321}^- + \bar C_{53} + \bar C_{42}\\
 \partial(C_{54}) & = & G_{1243} + G_{2143} + H_4 + F_4 + \bar C_{54} + \bar C_{43} 
\end{eqnarray*}
\begin{eqnarray*}
 \partial(D_{12}^+) & = & J_{213}^- + K_3^+ + H_1 + H_2 + \bar D_{12}^+ + \bar D_{14}^+ \\
 \partial(D_{12}^-) & = & I_1 + J_{213}^+ + K_3^- + H_1 + H_2 + \bar D_{12}^- + \bar D_{14}^- \\
 \partial(D_{13}^+) & = & J_{123}^- + J_{213}^+ + J_{312}^- + H_1 + H_3 + \bar D_{13}^+ + \bar D_{24}^+ \\
\partial(D_{13}^-) & = & J_{123}^+ + J_{213}^- + J_{312}^+ + H_1 + H_3 + \bar D_{13}^- + \bar D_{24}^- \\
 \partial(D_{14}^+) & = & J_{132}^- + K_2^+ + J_{312}^+ + H_1 + H_4 + \bar D_{14}^+ + \bar D_{34}^+ \\
 \partial(D_{14}^-) & = & J_{132}^+ + K_2^- + J_{312}^- + H_1 + H_4 + \bar D_{14}^- + \bar D_{34}^- \\
 \partial(D_{23}^+) & = & J_{123}^+ + J_{321}^- + H_2 + H_3 + \bar D_{23}^+ + \bar D_{12}^+ \\
 \partial(D_{23}^-) & = & J_{123}^- + I_2 + J_{321}^+ + H_2 + H_3 + \bar D_{23}^- + \bar D_{12}^- \\
 \partial(D_{24}^+) & = & J_{132}^+ + J_{231}^- + J_{321}^+ + H_2 + H_4 + \bar D_{24}^+ + \bar D_{13}^+ \\
 \partial(D_{24}^-) & = & J_{132}^- + J_{231}^+ + J_{321}^- + H_2 + H_4 + \bar D_{24}^- + \bar D_{13}^- \\
 \partial(D_{34}^+) & = & K_1^+ + J_{231}^+ + H_3 + H_4 + \bar D_{34}^+ + \bar D_{23}^+ \\
 \partial(D_{34}^-) & = & K_1^- + J_{231}^- + H_3 + H_4 + I_3 + \bar D_{34}^- + \bar D_{23}^- \\
 \partial(E_+) & = & K_1^+ + K_2^+ + K_3^+ \\
 \partial (E_1) & = & K_1^+ + K_2^- + K_3^- + W_1^+ + \bar E_1 + \bar E_3\\
 \partial(E_2) & = & K_1^- + K_2^+ + K_3^- + W_1^- + W_2^+ + \bar E_2 + \bar E_1\\
 \partial(E_3) & = & K_1^- + K_2^- + K_3^+ + W_2^- + \bar E_3 + \bar E_2 \\
& & \\
\partial(\bar A_{23}) & = & \bar C_{12} + \bar C_{21} + \bar C_{45} + \bar C_{54} + \bar D_{24}^- \\
 \partial(\bar A_{24}) & = & \bar C_{13} + \bar C_{53} + \bar D_{23}^+ + \bar D_{14}^- + \bar D_{24}^+ + \bar D_{34}^- + \bar E_3\\
 \partial(\bar A_{25}) & = & \bar C_{14} + \bar C_{35} + \bar D_{23}^- + \bar D_{24}^+ + \bar D_{14}^- + \bar D_{34}^+ + \bar E_1 \\
 \partial(\bar A_{26}) & = & \bar C_{15} + \bar C_{34} + \bar C_{43} + \bar D_{24}^- + \bar C_{51} \\
 \partial(\bar A_{34}) & = & \bar C_{31} + \bar C_{52} + \bar D_{13}^+ + \bar D_{14}^+ + \bar D_{23}^- + \bar D_{34}^- + \bar E_2 \\
 \partial(\bar A_{35}) & = & \bar D_{12}^+ + \bar D_{13}^- + \bar D_{23}^+ \\
 \partial(\bar A_{36}) & = & \bar C_{42} + \bar D_{12}^+ + \bar D_{14}^- + \bar D_{24}^+ + \bar D_{34}^- + \bar E_3 + \bar C_{52} \\
 \partial(\bar A_{45}) & = & \bar C_{25} + \bar C_{41} + \bar D_{12}^- + \bar D_{14}^+ + \bar D_{34}^- + \bar E_2 + \bar D_{24}^+ \\
 \partial(\bar A_{46}) & = & \bar C_{24} + \bar D_{12}^- + \bar D_{13}^+ + \bar D_{14}^- + \bar D_{34}^+ + \bar E_1 + \bar C_{53} \\
 \partial(\bar A_{56}) & = & \bar C_{23} + \bar C_{32} + \bar C_{51} + \bar D_{13}^- + \bar C_{54}
\end{eqnarray*}
\end{proposition}
The matrix of this operator is given in pages \pageref{page5to4}--\pageref{page5to4end}.
\hfill $\Box$

\begin{corollary}
\label{corhom5}
The kernel of the operator $\partial_5: C_5 \to C_4$ is generated by right parts of the formulas given in the following Proposition \ref{6to5}. In particular, $H_5(\overline{TD}_2(S^1), \Z_2) \simeq 0$. \hfill $\Box$
\end{corollary}

\begin{proposition}
\label{6to5}
The boundary operator $\partial_6: C_6 \to C_5$ of the mod 2 cellular complex of the CW-complex $\overline{TD}_2(S^1)$ is defined by following formulas:
\begin{eqnarray*}
\label{Afirst} \partial(A_{23}) & = & C_{12} + C_{21} + B + C_{45} + C_{54} + \bar A_{23} + \bar A_{26} \\
 \partial(A_{24}) & = & C_{13} + B + C_{53} + D_{23}^+ + \bar A_{24} + \bar A_{36} \\
 \partial(A_{25}) & = & C_{14} + B + C_{35} + D_{23}^- + D_{24}^+ + \bar A_{25} + \bar A_{46} \\
 \partial(A_{26}) & = & C_{15} + C_{34} + C_{43} + D_{24}^- + \bar A_{26} + \bar A_{56} \\
 \partial(A_{34}) & = & B + C_{31} + C_{52} + D_{13}^+ + D_{23}^- + \bar A_{34} + \bar A_{45} \\
 \partial(A_{35}) & = & B + D_{12}^+ + D_{13}^- + D_{14}^+ + D_{23}^+ + D_{24}^- + D_{34}^+ + E_+ \\
 \partial(A_{36}) & = & C_{42} + D_{12}^+ + D_{14}^- + D_{24}^+ + D_{34}^- + E_3 + \bar A_{36} + \bar A_{34} \\
 \partial(A_{45}) & = & B + C_{25} + C_{41} + D_{12}^- + D_{14}^+ + D_{34}^- + E_2 + \bar A_{45} + \bar A_{25} \\
 \partial(A_{46}) & = & C_{24} + D_{12}^- + D_{13}^+ + D_{14}^- + D_{34}^+ + E_1 + \bar A_{46} + \bar A_{24} \\
\label{Alast} \partial(A_{56}) & = & C_{23} + C_{32} + C_{51} + D_{13}^- + \bar A_{56} + \bar A_{23}
\end{eqnarray*}
The matrix of this operator is given in page \pageref{page6to5}. \hfill $\Box$
\end{proposition}

\begin{corollary}
\label{corhom6}
$H_6(\overline{TD}_2(S^1), \Z_2) \simeq 0$. \hfill $\Box$
\end{corollary}

\begin{corollary}
\label{corhom4}
$H_4(\overline{TD}_2(S^1), \Z_2) \simeq 0$.
\end{corollary}

\noindent
{\it Proof.} The Euler characteristic of the complex $\overline{TD}_2(S^1)$ is zero.
By Corollaries \ref{corhom1}, \ref{corhom2}, \ref{corhom3}, \ref{corhom5}, and \ref{corhom6}, the dimensions of the mod 2 homology groups in dimensions 0, 1, 2, 3, 5, and 6 satisfy the relation $\beta_0-\beta_1+\beta_2-\beta_3 -\beta_5+\beta_6 =0$. Thus, also $\beta_4 =0$. \hfill $\Box$
\medskip

Theorem \ref{mainthm1} follows immediately from Corollaries \ref{corhom6}, \ref{corhom5}, \ref{corhom4}, \ref{corhom3}, \ref{corhom2}, and \ref{corhom1}. \hfill $\Box$

\begin{proposition}
\label{pro17}
The subcomplex of \ $\overline{TD}_2(S^1)$ consisting of algebras of multiplicity five has only two non-trivial homology groups with coefficients in $\Z_2$: it are the groups $H_0$ and $H_1$, which are both isomorphic to $ \Z_2$. The one-dimensional homology group
of this subcomplex is generated by either the cell $\bar U_2$, or the cell $\bar Y_1$, or the cell $\nabla$.
\end{proposition}

\noindent
{\it Proof} follows immediately from the formulas of Propositions \ref{2to1}, \ref{3to2}, \ref{4to3}, and \ref{5to4}. \hfill $\Box$

\begin{corollary}
\label{cor99}
The homology groups mod 2 of the quotient complex of \ $\overline{TD}_2(S^1)$ \ by the subcomplex of algebras of multiplicity five are as follows: 
$H_0 \simeq 0, $ $H_1 \simeq \Z_2$, $H_2 \simeq \Z_2^3$, $H_3 \simeq \Z_2^2$. 
\end{corollary}

\noindent
{\it Proof.} It follows from Propositions \ref{2to1} and \ref{pro17} that the generator of the group $H_1$ of the subcomplex is the boundary of a two-dimensional chain in $\overline{TD}_2(S^1)$. The statement of Corollary \ref{cor99} follows then from the exact sequence of the pair. \hfill $\Box$
}

\section{Proof of Theorem \ref{mainthm2}}
\label{pro2}

\begin{proposition}
\label{W1}
The first Stiefel--Whitney class $w_1(\N_2)$ of the canonical normal bundle on $\overline{TD}_2(S^1)$ takes the non-zero value on the homology class of the cycle $\bar \Theta$.

The classes $w_2(\N_2)$ and $w_1^2(\N_2)$ take the non-zero value on the homology class of the cycle $\bar \Omega$.
\end{proposition}

\noindent
{\it Proof.} The closure of the cell $\bar \Omega$ is diffeomorphic to $\RP^2$, see \S~\ref{1sup}. It consists of the two-dimensional cell $\bar \Omega$, the one-dimensional cell $\bar \Theta$, and the $0$-dimensional cell $\bar \nabla$.

In the restriction to this closure, the bundle $\N_2$ splits into the sum of a trivial two-dimensional bundle, which is normal to the space of functions with $f'(\bullet) = f''(\bullet)=0$, and the two-dimensional normal bundle of the tautological bundle on $\RP^2$. By the Whitney multiplication formula (see \cite{Mi}, \S~4), the total Stiefel--Whitney class of the bundle $\N_2$ in the restriction to this projective plane is equal to $(1+\tau)^{-1} = 1+\tau + \tau^2$, where $\tau$ is the generator of the one-dimensional cohomology group of $\RP^2$.
\hfill $\Box$
\medskip

Define the two-dimensional cycle $\tilde X \subset \overline{TD}_2(S^1)$ as the closure of the union of all algebras $X(\alpha; \varphi, \varphi+\pi)$, see page \pageref{xx}. It is fibered over the space $\RP^1$ of pairs of opposite points $(\varphi, \varphi+\pi) \subset S^1$. The fiber over such a pair is homeomorphic to $\RP^1$ and consists of all algebras defined by the conditions $f(\varphi) =f(\varphi+\pi), f'(\varphi)=0=f'(\varphi+\pi)$, and $\xi f''(\varphi) + \eta f''(\varphi+\pi)$. The algebras corresponding to the fractions $(\xi:\eta)  \in \RP^1$ with $\xi \neq 0 \neq \eta$ are of class $X$, and the algebras corresponding to $\xi = 0$ or $\eta=0$ are of class $Y$. 
This fiber bundle is non-orientable, so the cycle $\tilde X$ is homeomorphic to the Klein bottle. 

\begin{proposition}
\label{W2}
The cycle $\tilde X$ is homologous in $\overline{TD}_2(S^1)$ to the cycle $\bar X^+ + \bar X^-$.
\end{proposition}

\noindent
{\it Proof.} Consider the closure of the union of all algebras $X(\alpha; a, b)$ in $\overline{TD}_2(S^1)$. It is a three-dimensional variety fibered over the space $\mbox{Sym}^2(S^1)$ of unordered pairs of points $(a, b) \subset S^1$ (with exceptional fibers over the degenerate pairs $(a, a)$). This space $\mbox{Sym}^2(S^1)$ is homeomorphic to the closed Moebius band. The cycle $\tilde X$ and the closure of the variety $\bar X^+ \cup \bar X^-$ are the preimages of some two homological loops in this Moebius band. Thus, the homology between these cycles is defined by the preimage of the two-dimensional chain defining the homology between these loops.
\hfill $\Box$

\begin{proposition}
\label{prop17}
1. The class $w_2(\N_2)$ takes zero value on the cycle $\tilde X$. 

2. The class $w_1^2(\N_2)$ takes the non-zero value on it.
\end{proposition}

{\it Proof.}  The restriction of the bundle $\N_2$ to the cycle $\tilde X$ is the sum of three subbundles: 
\begin{itemize}
\item[(a)]
a line bundle that is normal to the bundle of hyperplanes in $C^\infty(S^1, \R)$ defined by the conditions $f(\varphi) = f(\varphi+\pi)$,
\item[(b)]
a two-dimensional bundle that is normal to the bundle of subspaces of codimension two defined by the conditions $f'(\varphi)=0$ and $f'(\varphi+\pi)=0$, 
\item[(c)]
a line bundle whose fiber over any point $X(\alpha; \varphi, \varphi+\pi) \in \tilde X$ is normal to the hyperplane defined by the condition $f''(\varphi) = \alpha f''(\varphi+\pi)$.
\end{itemize}

The bundles (a) and (b) do not depend on $\alpha$ and are lifted from some bundles on the base $\RP^1 = \{\varphi, \varphi+\pi\}$ of the cycle $\tilde X$. The first of these bundles on $\RP^1$ is nontrivial, and the second bundle is the sum of a trivial and a non-trivial one-dimensional bundles. Therefore, the total Stiefel--Whitney class of the restriction of $\N_2$ to $\tilde X$ is equal to that of the one-dimensional bundle (c). This implies immediately statement 1 of Proposition \ref{prop17}.

The first Stiefel-Whitney class of bundle (c) takes the non-zero value on every fiber of the fiber bundle $\tilde X \to \RP^1$, and  also on the cross-section 
$\{X(1; \varphi, \varphi+\pi)\}$ of this fiber bundle consisting of all algebras $X(\alpha; \varphi, \varphi+\pi)$ with $\alpha \equiv 1$. Indeed, the section of the vector bundle c) defined by the cosets of the function $\cos \varphi$ intersects its zero section at only one point $X(-1; \varphi, \varphi+\pi)$ of the fiber over any point $(\varphi, \varphi+\pi) \in \RP^1$, $\varphi \neq \pm \pi/2$, and also at only one point $X(1; \pi/2, 3\pi/2)$ of the cycle $\{X(1; \varphi, \varphi+\pi)\}$.  Therefore, in the restriction to $\tilde X$, the first Stiefel-Whitney class of the bundle (c) is equal to the intersection index with the cycle $\{X(1; \varphi, \varphi+\pi)\}$. The self-intersection index of this cycle in $\tilde X$ is non-zero, hence $\langle w_1^2(\N_2), [\tilde X] \rangle =1$.
\hfill $\Box$

\begin{proposition}
\label{zero42}
The characteristic classes $w_1^3(\N_2), $ $w_1(\N_2) \smile w_2(\N_2)$, and $w_3(\N_2)$ take the zero value on the cycle $(\ref{kr26})$.
\end{proposition}

\noindent
{\it Proof.} The cycle $(\ref{kr26}) = L_{123} + L_{231} + L_{312} + V_1^+ + V_2^+$ consists of two parts. The first summand, $L_{123} + L_{231} + L_{312},$ is the union of all trinitary algebras characterized by three ordered points $a, b, c \in S^1$, whose cyclical order in $S^1$ is compatible with the orientation of $S^1$; the algebra corresponding to such a triple consists of functions $f$ such that
\begin{equation}
f'(a)=f''(a)=0, \qquad f'(b)=0, \qquad f(c)=f(b),
\label{alga1}
\end{equation} 
see Fig.~\ref{p42} (left). The second summand, $V_1^+ + V_2^+$, is the union of algebras characterized by pairs of points $a, b \in S^1$ and coefficients $(\xi:\eta) \in \RP^1$. Namely, the algebra corresponding to such a set of data consists of functions $f$ such that \begin{equation}
\label{alga2}
f'(a) = f''(a)=0, \qquad f(b)=f(a), \qquad \xi f'''(a) + \eta f'(b) =0.
\end{equation}

{
\begin{figure}
\begin{picture}(10,10)
\put(5,5){\circle{10}}
\put(4.8,0.5){\makebox(0,0)[cc]{ $\ast$}}
\put(4.8,1.3){\makebox(0,0)[cc]{ $\ast$}}
\put(4.5,0){\makebox(0,0)[cc]{\tiny $a$}}
\put(8.5,9){\makebox(0,0)[cc]{\tiny $b$}}
\put(0,5.5){\makebox(0,0)[cc]{\tiny $c$}}
\put(7,8.3){\makebox(0,0)[cc]{ $\ast$}}
\put(7.2,8.3){\line(-2,-1){6.3}}
\end{picture} \qquad 
\begin{picture}(3,10)
\put(0,4){\Large $+$}
\end{picture}
\qquad 
\begin{picture}(15,10)
\put(5,5){\circle{10}}
\put(4.8,0.5){\makebox(0,0)[cc]{ $\ast$}}
\put(4.5,0){\makebox(0,0)[cc]{\tiny $a$}}
\put(2.3,9.5){\makebox(0,0)[cc]{\tiny $b$}}
\put(4.8,1.3){\makebox(0,0)[cc]{ $\ast$}}
\put(5.3,0.9){\line(-1,4){1.97}}
\put(4.7,0.9){\line(-1,4){1.9}}
\end{picture}
\caption{The cycle (\ref{kr26})}
\label{p42}
\end{figure}
}

This cycle (\ref{kr26}) is the space of a fiber bundle over the circle: the projection map sends any algebra (\ref{alga1}) or (\ref{alga2}) to the point $a$. 

 The rotation group of the circle acts freely on this fiber bundle and defines a trivialization of it. The vector bundle $\N_2$ over the cycle (\ref{kr26}) is induced by this trivialization from the restriction of this vector bundle to an arbitrary single fiber. This fiber is two-dimensional, therefore all characteristic cohomology classes of dimension three induced from it are equal to zero. \hfill $\Box$ 

\begin{proposition}
\label{zero66}
The characteristic classes $w_1^3(\N_2), $ $w_1(\N_2) \smile w_2(\N_2)$, and $w_3(\N_2)$ take zero value on the cycle $(\ref{kr66})+(\ref{kr67}).$
\end{proposition}

\begin{figure}
\begin{picture}(10,10)
\put(5,5){\circle{10}}
\put(9,5){\circle*{0.9}}
\put(9.5,5){\makebox(0,0)[cc]{ $\ast$}}
\put(1.7,9.5){\makebox(0,0)[cc]{\small $a$}}
\put(1.7,0.5){\makebox(0,0)[cc]{\small $b$}}
\put(9,5){\line(-1,0){4}}
\put(5,5){\line(-2,3){2.3}}
\put(5,5){\line(-2,-3){2.3}}
\bezier{300}(9,5)(5.85,6.72)(2.7,8.45)
\end{picture} \qquad 
\begin{picture}(3,10)
\put(0,4){\Large $+$}
\end{picture}
\qquad 
\begin{picture}(15,10)
\put(5,5){\circle{10}}
\put(9,5){\circle*{0.9}}
\put(9.5,5){\makebox(0,0)[cc]{ $\ast$}}
\put(1.7,9.5){\makebox(0,0)[cc]{\small $b$}}
\put(1.7,0.5){\makebox(0,0)[cc]{\small $a$}}
\put(9,5){\line(-1,0){4}}
\put(5,5){\line(-2,3){2.3}}
\put(5,5){\line(-2,-3){2.3}}
\bezier{300}(9,5)(5.85,3.28)(2.7,1.55)
\end{picture}
\caption{The cycle (\ref{kr66})+(\ref{kr67})}
\label{p82}
\end{figure}

\noindent
{\it Proof.} This cycle
\begin{equation}
\label{eqJ}
\bar J_{312}^+ + \bar J_{312}^- + \bar J_{321}^+ + \bar J_{321}^-
\end{equation}
is the closure of the set of trinitary algebras characterized by pairs of distinct points $a, b \in (0,2\pi)$ and coefficients $(\sigma:\tau) \in \RP^1$: the algebra corresponding to these data is distinguished by the conditions 
\begin{equation}
\label{regJ}
f'(0)=0, \qquad f(a) = f(0) = f(b), \qquad \sigma f''(0) + \tau f'(a) = 0, 
\end{equation}
see Fig.~\ref{p82}. Denote this closure by  
\fbox{$\boxtimes$}.

Consider the unique continuous map \begin{equation}
\Pi: \fbox{$\boxtimes$} \to T^2
\end{equation}
 which sends any algebra of the form (\ref{regJ}) to the point $(a, b) \in (0,2\pi) \times (0, 2\pi)$.

\begin{lemma}
\label{lem44}
This map $\Pi$ is a locally trivial fiber bundle with the base $T^2$ and  the fibers homeomorphic to $S^1$. The closure of the set of all algebras $(\ref{regJ})$ with $\tau=0$ is a cross-section of this fiber bundle. 
\end{lemma}

\begin{figure}
\unitlength 1.1 mm
\begin{picture}(51,51)
\put(2,2){\vector(1,0){49}}
\put(2,2){\vector(0,1){49}}
\put(2,2){\line(1,1){40}}
\put(42,2){\line(0,1){40}}
\put(2,42){\line(1,0){40}}
\put(0.5,0.5){\makebox(0,0)[cc]{\small $\Theta$}}
\put(49,0){\makebox(0,0)[cc]{\small $a$}}
\put(0.5,49){\makebox(0,0)[cc]{\small $b$}}
\put(42,0){\makebox(0,0)[cc]{\small $2\pi$}}
\put(0,42){\makebox(0,0)[cc]{\small $2\pi$}}
\put(23,-0.6){\makebox(0,0)[cc]{\small $\bar V_2^{\pm}$}}
\put(23,44.3){\makebox(0,0)[cc]{\small $\bar V_2^{\pm}$}} 
\put(23.5,20){\makebox(0,0)[cc]{\small $\bar X^{\pm}$}}
\put(0,20){\makebox(0,0)[cc]{\small $\bar S_2$}}
\put(44,20){\makebox(0,0)[cc]{\small $\bar S_2$}}
\put(12,32){\makebox(0,0)[cc]{\small $\bar J^{\pm}_{321}$}}
\put(32,12){\makebox(0,0)[cc]{\small $\bar J^{\pm}_{312}$}}
\end{picture}
\caption{Base of the cycle (\ref{kr66})+(\ref{kr67})}
\label{base}
\end{figure}

\noindent
{\it Proof.} The target space $T^2$ of this map consists of six open strata, see Fig.~\ref{base}: 
\begin{itemize}
\item
two triangles distinguished by the inequalities $0 < a < b <2\pi$ and $0<b<a< 2\pi$, 
\item
three intervals 
$\{a=0, 0<b<2\pi\}$, $\{b=0, 0<a<2\pi\}$, and $\{0 < a = b <2\pi\}$,
\item 
the point $a = b = 0 \equiv 2\pi$.
\end{itemize}
Let us consider the preimages of all these strata.
\smallskip

\noindent
{\bf J}. The preimage $\Pi^{-1}(a,b)$ of any point $(a \neq b)$ of a two-dimensional stratum is the set of all algebras (\ref{regJ}) with arbitrary $(\sigma: \tau) \in \RP^1$. It splits into 
two intervals consisting of algebras from the cells $\bar J_{312}^\pm$ or $\bar J_{321}^\pm$ and two points of classes $\bar M$ and $\bar N$ corresponding to the values 
$(\sigma:\tau) =(0: 1)$ and $(1: 0)$. 
\medskip

\noindent
{\bf S}. The fiber $\Pi^{-1}(0, b)$ over any point $( 0, b)$, $ b \in (0, 2\pi)$, is the union of algebras satisfying the conditions
\begin{equation}
\label{stratS}
f'(0)=f''(0)=0, \qquad f(b)=f(0), \qquad \xi f^{\mbox{\i v}}(0) + \eta f'''(0)=0
\end{equation}
for some $(\xi:\eta) \in \RP^1$. These algebras belong to the cell $\bar S_2$ if $\xi \neq 0$ and to the cell $\bar U_2$ if $\xi=0$. Let us describe a local trivialization of the map $\Pi$ in the preimage of a neighborhood of the point $(0,b)$, that is, a homeomorphism between this preimage and the direct product of the fiber and this neighborhood. Any such trivialization can be defined by the system of its {\em leaves}, that is, the preimages of different points of the fiber under the natural projection of this direct product structure. So let us describe, for any $a \approx 0$, $\tilde b \approx b$, and $(\xi:\eta) \in \RP^1$, the intersection point of the fiber $\Pi^{-1}(a, \tilde b)$ with the leaf of our trivialization that contains the algebra (\ref{stratS}). These intersection points with the fibers $\Pi^{-1}(0, \tilde b),$ $\tilde b \approx b$, over the same stratum of $T^2$ are the similar algebras defined by the conditions
$$f'(0)=f''(0)=0, \qquad f(\tilde b)=f(0), \qquad \xi f^{\mbox{\i v}}(0) + \eta f'''(0)=0$$
with the same $(\xi:\eta)$.
Further, the intersection of such a leaf with the fiber $\Pi^{-1}(a, \tilde b)$, where $a \approx 0$, $\tilde b \approx b$, but $a \neq  0$, can be chosen as the algebra defined by the conditions
\begin{equation}
\label{leafS}
f'(0)=0, \qquad f(a)=f(0) =f(\tilde b), \qquad \left(4 \xi a- \eta a^2\right) f''(0) + 8 \xi f'(a)=0.
\end{equation}
Indeed, for $\xi \neq 0$ we take $\theta = - \frac{\eta}{\xi}$ and see from \S~\ref{incJS} that the algebra (\ref{stratS}) is the limit of algebras (\ref{eqJJ}), $a \to 0$, with $\lambda(a) = - \frac{1}{2}a -\frac{1}{8}\theta a^2.$ For $\xi =0$ the formula (\ref{leafS}) holds by the continuity.
\medskip

\noindent
{\bf X.} The fiber $\Pi^{-1}(a, a),$ $a\in (0,2\pi)$, consists of all algebras satisfying the conditions 
\begin{equation}
f'(0)=f'(a)=0, \qquad f(a)=f(0), \qquad \xi f''(a) + \eta f''(0) =0
\label{stratX}
\end{equation}
with arbitrary $(\xi: \eta) \in \RP^1$. These algebras are of class $\bar X^{\pm}$ if $\xi \neq 0 \neq \eta$, of class $\bar Y_2$ if $\xi=0$, and of class $\bar Y_1$ if $\eta=0$.
For any such algebra, the containing it leaf of the desired trivialization of the fiber bundle $\Pi$ intersects the fiber $\Pi^{-1}(\tilde a, \tilde a)$ with $\tilde a \approx a$ at the algebra defined by the similar condition 
$$f'(0)=f'(\tilde a)=0, \qquad f(\tilde a)=f(0), \qquad \xi f''(\tilde a) + \eta f''(0) =0$$
with the same $(\xi:\eta) \in \RP^1$. 
For any point $(\tilde a, \tilde b)$ with $\tilde a \approx a \approx \tilde b$, $\tilde a \neq \tilde b,$ the intersection of this leaf with the fiber $\Pi^{-1}(\tilde a, \tilde b)$ is the algebra defined by the conditions
\begin{equation}
\label{leafX}
f'(0)=0, \qquad f(\tilde a)=f(\tilde b)=f(0), \qquad  (\tilde b - \tilde a) \eta f''(0) - \xi f'(\tilde a)=0.
\end{equation}

\noindent
{\bf V.} For any $a \in (0, 2\pi)$, the fiber $\Pi^{-1}(a,0)$ is the union of algebras defined by the conditions
\begin{equation}
\label{stratV}
f'(0)= f''(0)=0, \qquad f(a) = f(0), \qquad \xi f'''(0)+\eta f'(a) =0 
\end{equation}
for some $(\xi:\eta) \in \RP^1$. These algebras are of classes $\bar V_2^{\pm}$ if $\xi  \neq 0 \neq \eta$, of class $\bar Y_2$ if $\xi=0$ and of class $\bar U_2$ if $\eta=0$. 
The desired leaf passing through the algebra (\ref{stratV}) intersects any neighboring fiber $\Pi^{-1}(\tilde a, 0)$, $\tilde a \approx a$, at the similar algebra 
$$f'(0)=0, \qquad f''(0)=0, \qquad f(\tilde a) = f(0), \qquad \xi f'''(0)+\eta f'(\tilde a) =0 $$
with the same $(\xi:\eta) \in \RP^1$.
For any neighboring point $(\tilde a, b)$ with $\tilde a \approx a$, $b \approx 0$, $b\neq 0$, this leaf intersects the fiber $\Pi^{-1}(\tilde a, b)$ at the point corresponding to the algebra defined by the conditions
\begin{equation}
\label{leafV}
f'(0)=0, \qquad f(a)=f(0)=f(b), \qquad 3 \xi f''(0) - b \eta f'(a) = 0.
\end{equation} 

\noindent
${\bf \Theta}.$
The fiber $\Pi^{-1}(0,0)$ consists of algebras defined by the conditions 
\begin{equation}
\label{stratTheta}
f'(0)=f''(0)=f'''(0)=0, \qquad \xi f^{\tiny \mbox{V}}(0) + \eta f^{\tiny \mbox{IV}}(0)=0
\end{equation}
with arbitrary $(\xi:\eta) \in \RP^1$. 
These algebras are of class $\Theta$ if $\xi=0$ or of class $\nabla$ if $\xi=0$. 
The leaf of the local trivialization of the bundle $\Pi$ that passes through the point (\ref{stratTheta}) intersects any fiber $\Pi^{-1}(a, b)$, where $0 \neq a \neq b \neq 0$, at the point (\ref{regJ}) with coefficients $(\sigma: \tau) $ satisfying the equality 
\begin{equation}
\label{leafTheta}
2\left(5 \xi b - \eta (ab+b^2)\right) \sigma + a(a-b) \left(5\xi- \eta (2a+b)\right)\tau =0.
\end{equation}
It is easy to check that for any $a$ and $b$ such that $0 \neq a \neq b \neq 0$ this formula defines a diffeomorphism between the fibers over the points $(0, 0)$ and $(a,b)$.
For the fibers over the points of the neighboring strata $\{a =0\},$ $\{b=0\}$ and $\{a=b\}$ these intersection points are defined by continuity.

So, the map $\Pi: \fbox{$\boxtimes$} \to T^2$ is indeed a locally trivial  fiber bundle. The  set of algebras $(\ref{regJ})$ with $\tau=0$ is a cross-section of this fiber bundle, as follows from formulas (\ref{leafS}), (\ref{leafX}), (\ref{leafV}), and (\ref{leafTheta}). Namely, this cross-section maps the points of the non-typical strata of $T^2$  to  the following points of fibers: $\xi=0$ in (\ref{stratS}), $\xi=0$ in (\ref{stratX}), $\eta=0$ in (\ref{stratV}), and $\xi=0$ in (\ref{stratTheta}). \hfill $\Box$

\begin{lemma}
\label{lem55}
The fiber bundle $\Pi: \fbox{$\boxtimes$} \to T^2$ considered in Lemma \ref{lem44} is non-orientable: it changes the orientation of fibers over a loop in the base if and only if this loop has an odd intersection index with the set of points $(0, b)$. 
\end{lemma}

\noindent
{\it Proof.} For any point $(a, b) \in T^2$ with $0 \neq a \neq b \neq 0$,
define the standard orientation of the fiber $\Pi^{-1}(a, b)$ by the increase of the parameter $\sigma/\tau$ in (\ref{regJ}). The formulas (\ref{leafS}), (\ref{leafX}), and (\ref{leafV}) imply that crossing the strata $\{b=0\}$ and $\{a=b\}$ of the base changes this orientation, and crossing the stratum $\{a=0\}$ preserves it. \hfill $\Box$

\begin{corollary}
\label{cor77}
There is a homeomorphism 
\begin{equation}
\fbox{$\boxtimes$} \simeq S^1 \times K^2
\label{propro}
\end{equation}
of the space \fbox{$\boxtimes$} to the product of the circle and the Klein bottle: the
composition of this homeomorphism and the projection of the product to the first factor $S^1$ maps any algebra defined by conditions $($\ref{regJ}$)$, $($\ref{stratS}$)$, $($\ref{stratX}$)$, $($\ref{stratV}$)$, or $($\ref{stratTheta}$)$ to the point $b$ participating in these conditions. \hfill $\Box$ 
\end{corollary}

The restriction of the bundle $\N_2$ to the variety \fbox{$\boxtimes$} splits into the sum of four line bundles. Namely, on the open subset in  \fbox{$\boxtimes$} consisting 
of the algebras (\ref{eqJ}) with $a \neq 0 \neq b$, these bundles are normal to the bundles of hyperplanes in $C^{\infty}(S^1, \R)$ defined by the four equations (\ref{eqJ}).  These bundles can extended by continuity to all remaining points of the variety \fbox{$\boxtimes$}. 

The first line bundle is trivial by definition. The second and the third  bundles are lifted from some two bundles on the base $T^2$ of the map $\Pi$.

\begin{lemma}
The first Stiefel--Whitney classes of the latter two one-dimen\-si\-onal bundles are trivial.
\end{lemma}

\noindent
{\it Proof.} The group $H_1(T^2, \Z_2)$ is generated by two loops: one of them consists of all points $(a, b_0)$ for some fixed generic $b_0 \in S^1$, and the other loop consists of all points $(a_0, b)$ for some $a_0 \in S^1$. 
Any function $f: S^1 \to \R$ defines the cross-sections of both these line bundles.
Suppose that $f(a_0) \neq f(0) \neq f(b_0)$. 
The value of the first Stiefel--Whitney class of the second (respectively, the third) line bundle on the first loop is equal to the parity of the number of points $a \in S^1$ (including $a=0$ and counted with multiplicities) such that $f(a)=f(0)$ (respectively, $f(a_0)=f(0)$). Both these numbers are even (and the second of them is just zero). The proof for the second loop is analogous. \hfill $\Box$

\begin{corollary}
\label{cor88}
The total Stiefel--Whitney class of the restriction of $N_2$ to the variety \fbox{$\boxtimes$} coincides with that of the fourth line bundle. In particular, the classes $w_2(N_2)$ and $w_3(N_2)$ are trivial on this variety. 
\hfill $\Box$
\end{corollary}

\begin{lemma}
\label{lem66}
1. The first Stiefel--Whitney class of the fourth line bundle takes the non-zero value on any fiber of the fiber bundle $\Pi : \fbox{$\boxtimes$} \to T^2  $. 

2. This class takes the zero value on any loop in the cross-section of this fiber bundle mentioned in  Lemma \ref{lem44}.
\end{lemma}

\noindent
{\it Proof.} 1. Let $(a,b)$ be an arbitrary point of $T^2$ such that $0 \neq a \neq b \neq 0$.
Let $f: S^1 \to \R$ be any function with $f''(0) \neq 0 \neq f'(a)$. The value of the first Stiefel--Whitney class of the fourth line bundle on the cycle $\Pi^{-1}(a,b)$ is equal to the parity of the number of points $(\sigma:\tau) \in \RP^1$ such that $\sigma f''(0) + \tau f'(a) =0$. This equality holds at the single point $(\sigma:\tau) \equiv (-f'(a):f''(0))$.

2. Any function $f$ with $f''(0) \neq 0$ defines an everywhere non-zero cross-section of this line bundle over the entire this cross-section of $\Pi$. \hfill $\Box$

\begin{corollary}
\label{cor22}
The characteristic class $w_1^2(N_2)$ takes the non-zero value on any fiber $K^2$  of the projection $\fbox{$\boxtimes$}   \to S^1$ mentioned in Corollary \ref{cor77}. 
The restriction of the class $w_1(N_2)$ to the variety \fbox{$\boxtimes$} is induced
by the projection of the  product $($\ref{propro}$)$ to the second factor
from the restriction of this class to this fiber.
In particular, the class $w_1^3(N_2)$ takes the zero value on the cycle \fbox{$\boxtimes$}. \hfill $\Box$
\end{corollary}

\noindent
{\it Proof.} These facts easily follow from Lemmas \ref{lem55} and \ref{lem66}. \hfill $\Box$
\medskip

Proposition \ref{zero66} follows immediately from  Corollaries \ref{cor88} and \ref{cor22}. \hfill $\Box$

\medskip
By Corollary \ref{corhom3}, the group $H_3(\overline{TD}_2(S^1), \Z_2)$ is generated by the cycles (\ref{kr26}) and (\ref{kr66})+(\ref{kr67}). Therefore, Theorem \ref{mainthm2} is covered by Propositions \ref{W1}, \ref{W2}, \ref{zero42}, and \ref{zero66}. \hfill $\Box$

\newpage
\section*{Appendix: matrices}

\FloatBarrier

\begin{table}
\caption{Boundary operator $\partial_2: C_2 \to C_1$}
\label{page2to1}
\begin{center}
\begin{tabular}{| l || c | c | c | c|c|c|c|}
\hline
$\vphantom{\overline{(}}$ & $\bar Y_1 \vphantom{\overline{(}}$ & $\bar Y_2$ & $\bar U_1$ & $\bar U_2$ & $\bar Z$ & $\bar \Theta$ & $\nabla$ \\
\hline 
\hline
$Y_1 $ & + & + & & & & & + \\
\hline
$Y_2$ & + & + & & & & & + \\
\hline
$U_1$ & & & + & + & & & + \\
\hline
$U_2$ & & & + & + & & & + \\
\hline
$Z$ & & & & & & & + \\
\hline
$\Theta$ & & & & & & & \\
\hline
$\bar L_{123} \vphantom{\overline{(}}$ & & + & + & & + & & \\
\hline
$\bar L_{231} \vphantom{\overline{(}}$ & + & & & + & + & & \\
\hline
$\bar L_{312} \vphantom{\overline{(}}$ & & + & & + & + & & \\
\hline
$\bar L_{132} \vphantom{\overline{(}}$ & + & & & + & + & & \\
\hline
$\bar L_{213} \vphantom{\overline{(}}$ & & + & + & & + & & \\
\hline
$\bar L_{321} \vphantom{\overline{(}}$ & & + & & + & + & & \\
\hline
$\bar M_1 \vphantom{\overline{(}}$ & + & & + & + & & & \\
\hline
$\bar M_2 \vphantom{\overline{(}}$ & + & & + & + & & & \\
\hline
$\bar M_3 \vphantom{\overline{(}}$ & & + & & & & & \\
\hline
$\bar N_1 \vphantom{\overline{(}}$ & + & + & & + & & & \\
\hline
$\bar N_2 \vphantom{\overline{(}}$ & + & + & & + & & & \\
\hline
$\bar N_3 \vphantom{\overline{(}}$ & & & + & & & & \\
\hline
$\bar P_1 \vphantom{\overline{(}}$ & & & & & & & \\
\hline
$\bar P_2 \vphantom{\overline{(}}$ & & & & & & & \\
\hline
$\bar S_1 \vphantom{\overline{(}}$ & & & & & & & \\
\hline
$\bar S_2 \vphantom{\overline{(}}$ & & & & & & & \\
\hline
$\bar X^+ \vphantom{\overline{(}}$ & + & + & & & & & \\
\hline
$\bar X^- \vphantom{\overline{(}}$ & + & + & & & & & \\
\hline
$\bar V_1^+ \vphantom{\overline{(}}$ & + & & + & & & & \\
\hline
$\bar V_2^+ \vphantom{\overline{(}}$ & & + & & + & & & \\
\hline
$\bar V_1^- \vphantom{\overline{(}}$ & + & & + & & & & \\
\hline
$\bar V_2^- \vphantom{\overline{(}}$ & & + & & + & & & \\
\hline
$\bar \Omega \vphantom{\overline{(}}$ & & & & & & & \\
\hline
\end{tabular}
\end{center}
\end{table}

{
\tiny

\noindent
\begin{table}
\caption{Boundary operator $\partial_3: C_3 \to C_2,$ left top part}
\label{3to2lt}
\begin{tabular}{| l || c | c | c | c|c|c|c|c|c|c|c|c|c|c|c|c|c|c|c|c|c|c|c|}
\hline
$\vphantom{\overline{(}} $ & $Y_{1}$ & $Y_2$ & $U_1$ & $U_2$ & $Z$ & $\Theta$ & $\bar L_{123}$ & $\bar L_{231}$ & $\bar L_{312}$ & $\bar L_{132}$ & $\bar L_{213}$ & $\bar L_{321}$ & $\bar M_1$ & $\bar M_2$ & $\bar M_3$ & $\bar N_1$ & $\bar N_2$ & $\bar N_3$ \\
\hline
$L_{123}$ & & & + & & + & & + & & + & & & & & & & & & \\
\hline
$L_{231}$ & + & & & + & & & + & + & & & & & & & & & & \\
\hline
$L_{312}$ & & + & & & + & & & + & + & & & & & & & & & \\
\hline
$L_{132}$ & + & & & & + & & & & & + & & + & & & & & & \\
\hline
$L_{213}$ & & + & + & & & & & & & + & + & & & & & & & \\
\hline
$L_{321}$ & & & & + & + & & & & & & + & + & & & & & & \\
\hline
$M_1$ & + & & + & & & & & & & & & & + & & + & & & \\
\hline
\end{tabular}
\end{table}
\label{page3to2beg}

\FloatBarrier

\newpage

{
\tiny

\begin{table}
\caption{Boundary operator $\partial_3: C_3 \to C_2,$ right top part}
\label{3to2rt}
\begin{tabular}{| l || c | c | c | c|c|c|c|c|c|c|c|c|c|c|c|c|c|c|c|c|c|c|c|}
\hline
$\vphantom{\overline{(}} $ & $\bar P_1$ & $\bar P_2$ & $\bar S_1$ & $\bar S_2$ & $\bar X^+$ & $\bar X^-$ & $\bar V_1^+$ & $\bar V_2^+$ & $\bar V_1^-$ & $\bar V_2^-$ & $\bar \Omega$ \\
\hline
$L_{123}$ & & & & & & & & & & & \\
\hline
$L_{231}$ & & & & & & & & & & & \\
\hline
$L_{312}$ & & & & & & & & & & & \\
\hline
$L_{132}$ & & & & & & & & & & & \\
\hline
$L_{213}$ & & & & & & & & & & & \\
\hline
$L_{321}$ & & & & & & & & & & & \\
\hline
$M_1$ & & & & & & & & & & & \\
\hline
\end{tabular}
\end{table}

\noindent
\begin{table}
\caption{Boundary operator $\partial_3: C_3 \to C_2,$ left bottom part}
\label{3to2lb}
\begin{tabular}{| l || c | c | c | c|c|c|c|c|c|c|c|c|c|c|c|c|c|c|c|c|c|c|c|}
\hline
$\vphantom{\overline{(}} $ & $Y_{1}$ & $Y_2$ & $U_1$ & $U_2$ & $Z$ & $\Theta$ & $\bar L_{123}$ & $\bar L_{231}$ & $\bar L_{312}$ & $\bar L_{132}$ & $\bar L_{213}$ & $\bar L_{321}$ & $\bar M_1$ & $\bar M_2$ & $\bar M_3$ & $\bar N_1$ & $\bar N_2$ & $\bar N_3$ \\
\hline
$\bar J_{123}^+$ & & & & & & & & & & & & & + & & & & & + \\
\hline
$\bar J_{123}^-$ & & & & & & & & & & & & & + & & & & & + \\
\hline
$\bar J_{132}^+$ & & & & & & & & & & & & & + & & & & + & \\
\hline
$\bar J_{132}^-$ & & & & & & & & & & & & & + & & & & + & \\
\hline
$\bar J_{213}^+$ & & & & & & & & & & & & & & + & & & & + \\
\hline
$\bar J_{213}^-$ & & & & & & & & & & & & & & + & & & & + \\
\hline
$\bar J_{231}^+$ & & & & & & & & & & & & & & + & & + & & \\
\hline
$\bar J_{231}^-$ & & & & & & & & & & & & & & + & & + & & \\
\hline
$\bar J_{312}^+$ & & & & & & & & & & & & & & & + & & + & \\
\hline
$\bar J_{312}^-$ & & & & & & & & & & & & & & & + & & + & \\
\hline
$\bar J_{321}^+$ & & & & & & & & & & & & & & & + & + & & \\
\hline
$\bar J_{321}^-$ & & & & & & & & & & & & & & & + & + & & \\
\hline
$\bar K_1^+$ & & & & & & & & & & & & & & & & & + & + \\
\hline
$\bar K_1^-$ & & & & & & & & & & & & & & & & & + & + \\
\hline
$\bar K_2^+$ & & & & & & & & & & & & & & & & + & & + \\
\hline
$\bar K_2^-$ & & & & & & & & & & & & & & & & + & & + \\
\hline
$\bar K_3^+$ & & & & & & & & & & & & & & & & + & + & \\
\hline
$\bar K_3^-$ & & & & & & & & & & & & & & & & + & + & \\
\hline
$\bar W_1^+$ & & & & & & & & & & & & & & & & & & \\
\hline
$\bar W_1^-$ & & & & & & & & & & & & & & & & & & \\
\hline
$\bar W_2^+$ & & & & & & & & & & & & & & & & & & \\
\hline
$\bar W_2^-$ & & & & & & & & & & & & & & & & & & \\
\hline
\end{tabular}
\end{table}
}

\bigskip

\noindent
\begin{table}
\caption{Boundary operator $\partial_3: C_3 \to C_2,$ right bottom part}
\label{3to2rb}
\begin{tabular}{| l || c | c | c | c|c|c|c|c|c|c|c|c|c|c|c|c|c|c|c|c|c|c|}
\hline
$\vphantom{\overline{(}} $ & $\bar P_1$ & $\bar P_2$ & $\bar S_1$ & $\bar S_2$ & $\bar X^+$ & $\bar X^-$ & $\bar V_1^+$ & $\bar V_2^+$ & $\bar V_1^-$ & $\bar V_2^-$ & $\bar \Omega$ \\
\hline
$\bar J_{123}^+$ & & & & & & + & & + & & & \\
\hline
$\bar J_{123}^-$ & & & + & & + & & & & & + & \\
\hline
$\bar J_{132}^+$ & & & & + & + & & & & + & & \\
\hline
$\bar J_{132}^-$ & & & & & & + & + & & & & \\
\hline
$\bar J_{213}^+$ & & & + & & + & & & & & + & \\
\hline
$\bar J_{213}^-$ & & & & & & + & & + & & & \\
\hline
$\bar J_{231}^+$ & & & & & & + & + & & & & \\
\hline
$\bar J_{231}^-$ & & & & + & + & & & & + & & \\
\hline
$\bar J_{312}^+$ & & & & & & + & & + & & & \\
\hline
$\bar J_{312}^-$ & & & & + & + & & & & & + & \\
\hline
$\bar J_{321}^+$ & & & & + & + & & & & & + & \\
\hline
$\bar J_{321}^-$ & & & & & & + & & + & & & \\
\hline
$\bar K_1^+$ & & & & & & & + & + & & & \\
\hline
$\bar K_1^-$ & & + & & & & & & & + & + & \\
\hline
$\bar K_2^+$ & & & & & & & + & + & & & \\
\hline
$\bar K_2^-$ & & + & & & & & & & + & + & \\
\hline
$\bar K_3^+$ & & & & & & & & & & & \\
\hline
$\bar K_3^-$ & + & & & & & & & & & & \\
\hline
$\bar W_1^+$ & + & & & & & & + & & + & & + \\
\hline
$\bar W_1^-$ & + & & & & & & + & & + & & + \\
\hline
$\bar W_2^+$ & & + & & & & & & + & & + & + \\
\hline
$\bar W_2^-$ & & + & & & & & & + & & + & + \\
\hline
\end{tabular}
\end{table}
\label{page3to2end}

\noindent

\newpage

{\small

\FloatBarrier

\begin{table}
\caption{Boundary operator $ \partial_5: C_5 \to C_4$, page 1}
\label{page5to4}
\begin{tabular}{| l || c | c | c | c|c|c|c|c|c|c|c|c|}
\hline
& $F_1$ & $F_2$ & $F_3$ & $F_4$ & $G_{1234}$ & $G_{1243}$ & $G_{1324}$ & $G_{1342}$ & $G_{1423}$ & $G_{1432}$ & $G_{2134}$ & $G_{2143}$ \\
\hline
$B$ & & & & & & & & & & & & \\
\hline
$C_{12}$ & + & & & & + & + & & & & & & \\
\hline
$C_{13}$ & & & & & & & & + & & & & \\
\hline
$C_{14}$ & & & & & & & + & & & & & \\
\hline
$C_{15}$ & & & & & & & & & + & + & & \\
\hline
$C_{21}$ & + & & & & & & & & & & + & + \\
\hline
$C_{23}$ & & + & & & & & & & && & \\
\hline
\end{tabular}
\end{table}

\FloatBarrier

\begin{table}
\caption{Boundary operator $ \partial_5: C_5 \to C_4$, page 2}
\label{5to22}
\begin{tabular}{| l | c | c | c | c|c|c|c|c|c|c|c|c|}
\hline
$\vphantom{\overline{(}}$ & $G_{2341}$ &$G_{2431}$ & $G_{3142}$ & $G_{3241}$ & $H_1$ & $H_2$ & $H_3$ & $H_4$ & $I_1$ & $I_2$ & $I_3$ & $\bar B$ \\
\hline
$B$ & & & & & + & + & + & + & & & & \\
\hline
$C_{12}$ & & & & & + & & & & & & & \\
\hline
$C_{13}$ & & & & & + & & & & & & & \\
\hline
$C_{14}$ & & & & & + & & & & & & & \\
\hline
$C_{15}$ & & & & & & & & & & & & \\
\hline
$C_{21}$ & & & & & & + & & & & & & \\
\hline
$C_{23}$ & + & & & & + & + & & & && & \\
\hline
\end{tabular}
\end{table}
\newpage 

\begin{table}
\caption{Boundary operator $ \partial_5: C_5 \to C_4$, page 3}
\label{5to43}
\begin{tabular}{| l | c | c | c | c|c|c|c|c|c|c|c|c|}
\hline
& $J_{123}^+$ & $J_{132}^+$ & $J_{213}^+$ & $J_{231}^+$ & $J_{312}^+$ & $J_{321}^+$ & $J_{123}^-$ & $J_{132}^-$ & $J_{213}^-$ & $J_{231}^-$ & $J_{312}^-$ & $J_{321}^-$ \\
\hline
$B$ & & & & & & & & & & & & \\
\hline
$C_{12}$ & & & & & & & & & & & & \\
\hline
$C_{13}$ & + & & & & & & & & & & & \\
\hline
$C_{14}$ & & + & & & & & + & & & & & \\
\hline
$C_{15}$ & & & & & & & & + & & & & \\
\hline
$C_{21}$ & & & & & & & & & & & & \\
\hline
$C_{23}$ & + & & & & & & & & && & \\
\hline
\end{tabular}
\end{table}

\newpage

\begin{table}
\caption{Boundary operator $ \partial_5: C_5 \to C_4$, page 4}
\label{5to44}
\begin{tabular}{| l | c | c | c | c|c|c|c|c|c|c|c|c|}
\hline
& $K_{1}^+$ & $K_{2}^+$ & $K_{3}^+$ & $K_{1}^-$ & $K_{2}^-$ & $K_{3}^-$ & $W_1^+$ & $W_2^+$ & $W_1^-$ & $W_2^-$ \\
\hline
$B$ & & & & & & & & & & \\
\hline
$C_{12}$ & & & & & & & & & & \\
\hline
$C_{13}$ & & & & & & & & & & \\
\hline
$C_{14}$ & & & & & & & & & & \\
\hline
$C_{15}$ & & & & & & & & & & \\
\hline
$C_{21}$ & & & & & & & & & & \\
\hline
$C_{23}$ & & & & & & & & & & \\
\hline
\end{tabular}
\end{table}

\newpage

\begin{table}
\caption{Boundary operator $ \partial_5: C_5 \to C_4$, page 5}
\label{5to45}
\begin{tabular}{| l | c | c | c | c|c|c|c|c|c|c|c|c|}
\hline
$\vphantom{\overline{(}}$ & $\bar C_{12}$ & $\bar C_{13}$ & $\bar C_{14}$ & $\bar C_{15}$ & $\bar C_{21}$ & $\bar C_{23}$ & $\bar C_{24}$ & $\bar C_{25}$ & $\bar C_{31}$ & $\bar C_{32}$ & $\bar C_{34}$ & $\bar C_{35}$ \\
\hline
$B$ & & & & & & & & & & & & \\
\hline
$C_{12}$ & + & & & & & & & & & & & \\
\hline
$C_{13}$ & & + & & & & & & & & & & \\
\hline
$C_{14}$ & & & + & & & & & & & & & \\
\hline
$C_{15}$ & & & & + & & & & & & & & \\
\hline
$C_{21}$ & & & & + & + & & & & & & & \\
\hline
$C_{23}$ & + & & & & & + & & & && & \\
\hline
\end{tabular}
\end{table}

\newpage

\begin{table}
\caption{Boundary operator $ \partial_5: C_5 \to C_4$, page 6}
\label{5to46}
\begin{tabular}{| l | c | c | c | c|c|c|c|c|c|c|c|c|}
\hline
$\vphantom{\overline{(}}$ & $\bar C_{41}$ & $\bar C_{42}$ & $\bar C_{43}$ & $\bar C_{45}$ & $\bar C_{51}$ & $\bar C_{52}$ & $\bar C_{53}$ & $\bar C_{54}$ & $\bar D^+_{12}$ & $\bar D^+_{13}$ & $\bar D^+_{14}$ & $\bar D^+_{23}$ \\
\hline
$B$ & & & & & & & & & & & & \\
\hline
$C_{12}$ & & & & & + & & & & & & & \\
\hline
$C_{13}$ & & & & & & + & & & & & & \\
\hline
$C_{14}$ & & & & & & & + & & & & & \\
\hline
$C_{15}$ & & & & & & & & + & & & & \\
\hline
$C_{21}$ & & & & & & & & & & & & \\
\hline
$C_{23}$ & & & & & & & & & && & \\
\hline
\end{tabular}
\end{table}

\newpage

\begin{table}
\caption{Boundary operator $ \partial_5: C_5 \to C_4$, page 7}
\label{page5to4end}
\begin{tabular}{| l | c | c | c | c|c|c|c|c|c|c|c|c|}
\hline
$\vphantom{\overline{(}}$ & $\bar D^+_{24}$ & $\bar D^+_{34}$ & $\bar D^-_{12}$ & $\bar D^-_{13}$ & $\bar D^-_{14}$ & $\bar D^-_{23}$ & $\bar D^-_{24}$ & $\bar D^-_{34}$ & $\bar E_{+}$ & $\bar E_{3}$ & $\bar E_{2}$ & $\bar E_{1}$\\
\hline
$B$ & & & & & & & & & & & & \\
\hline
$C_{12}$ & & & & & & & & & & & & \\
\hline
$C_{13}$ & & & & & & & & & & & & \\
\hline
$C_{14}$ & & & & & & & & & & & & \\
\hline
$C_{15}$ & & & & & & & & & & & & \\
\hline
$C_{21}$ & & & & & & & & & & & & \\
\hline
$C_{23}$ & & & & & & & & & && & \\
\hline
\end{tabular}
\end{table}

\newpage 

\begin{table}
\caption{Boundary operator $ \partial_6: C_6 \to C_5$}
\label{page6to5}
\begin{tabular}{| l | c | c | c | c|c|c|c|c|c|c|}
\hline
& $A_{23}$ & $A_{24}$ & $A_{25}$ & $A_{26}$ & $A_{34}$ & $A_{35}$ & $A_{36}$ & $A_{45}$ & $A_{46}$ & $A_{56}$ \\
\hline
$B$ & + & + & + & & + & + & & + & & \\
\hline
$C_{12}$ & + & & & & & & & & & \\
\hline
$C_{13}$ & & + & & & & & & & & \\
\hline
$C_{14}$ & & & + & & & & & & & \\
\hline
$C_{15}$ & & & & + & & & & & & \\
\hline
$C_{21}$ & + & & & & & & & & & \\
\hline
$C_{23}$ & & & & & & & & & & + \\
\hline
\end{tabular}
\end{table}

}

\FloatBarrier

}
}

\end{document}